\documentclass[11pt,twoside]{article}
\usepackage{amsmath,amsfonts,amssymb,cases,mathdots,color,colordvi,url,tikz}
\usepackage{hyperref}
\usepackage{bm}
\usepackage{bbm}
\usepackage{pgfpages}
\usepackage{graphics}%\textwidth 16.cm
\usepackage{epsfig}
\usepackage{algorithm}
\usepackage{algorithmic}
\usepackage{subcaption,caption}
\usepackage{cleveref}
\newtheorem{proposition}{Proposition}%[section]
\newtheorem{theorem}[proposition]{Theorem}

\newtheorem{definition}[proposition]{Definition}
\newtheorem{example}[proposition]{Example}

\newtheorem{remark}{Remark}

\newcommand{\be}{\begin{equation}}
\newcommand{\ee}{\end{equation}}
\newcommand{\ba}{\begin{eqnarray}}
\newcommand{\ea}{\end{eqnarray}}
\newcommand{\bas}{\begin{eqnarray*}}
	\newcommand{\eas}{\end{eqnarray*}}

\graphicspath{ {./Figures/} }
\def\bfp{{\bf p}}

\def\bfr{{\bf r}}
\def\bfs{{\bf s}}

\def\bfd{{\bf d}}

\def\bfr{{\bf r}}

\def\bfx{{\bf x}}

\def\bfw{{\bf w}}

\def\bfone{{\bf 1}}

\def\diag{\mbox{diag}}
\def\Jen1{{J_{{\bf e}_{n+1} }}}
\def\JTen1{{J^\top_{{\bf e}_{n+1}}}}
\def\mdrp{\mbox{mdrp}}
\def\MDRP{\mbox{MDRP}}

\def\mvp{\mbox{mvp}}

\def\ef{\mbox{ef}}
\def\efdr{\mbox{dr}}

\def\cml{\mbox{cml}}
\def\ew{\mbox{ew}}
\def\mdp{\mbox{mdp}}
\def\MDP{\mbox{MDP}}
\def\mv{\mbox{mv}}

\title{An Optimization Study of Diversification Return Portfolios}

\author{Chao Ding\thanks{Institute of Applied Mathematics, Chinese Academy of Sciences, Beijing, China.
	Email: dingchao@amss.ac.cn.}
	\ \ \mbox{and} \ \
	Hou-Duo Qi\thanks{Department of Applied Mathematics, The Hong Kong Polytechnic University, Hong Kong. 
		Email: houduo.qi@polyu.edu.hk.}}

\date{\today}

\begin{document}
	
	\maketitle

\thispagestyle{empty}

\begin{abstract}
The concept of Diversification Return (DR) was introduced by Booth and Fama in 1990s and
it has been well studied in the finance literature mainly focusing on
the various sources it may be generated. 
However, unlike the classical Mean-Variance (MV) model of Markowitz,
DR portfolios lack optimization theory for justifying their often outstanding empirical performance.
In this paper, we first explain what the DR criterion tries to achieve in terms of portfolio centrality.
A consequence of this explanation is that practically imposed norm constraints in fact implicitly enforce constraints on DR.
%and this in turn justifies the use of norm constraints  in portfolio construction.
We then derive the maximum DR portfolio under given risk and obtain
the efficient DR frontier. 
We further develop a separation theorem for this frontier and
establish a relationship between the DR frontier  and Markowitz MV 
efficient frontier. 
In the particular case where the variance vector is proportional to the expected return vector of the underlining assets, the two frontiers yield same efficient portfolios.
The proof techniques heavily depend on recently developed geometric interpretation of the maximum DR portfolio. 
Finally, we use DAX30 stock data to illustrate the obtained results and demonstrate an interesting link to the maximum diversification ratio portfolio studied by Choueifaty and Coignard.

\vspace{3mm}

\vspace{3mm}
%\noindent{\bf \textbf{Mathematical Subject Classification}:} 90C26, 90C30, 90C90
\end{abstract}
{}
%\tableofcontents
%\numberwithin{equation}{section}

%%%%%%%%%%%%%%%%%%%%%%%%%%%%%%%%%%%%%%%%%%%%%%
\section{Introduction} \label{Section-Introduction}

Diversification Return (DR) of a portfolio was initially studied by Booth and Fama \cite{booth1992diversification}. 
The main purpose was to understand how much an excess rate of return could be
achieved by a portfolio when compared with the simply weighted return of
its constituents. Booth and Fama \cite{booth1992diversification} derived the DR
using the compound return, while Willenbrock \cite{willenbrock2011diversification} used the geometric return. 
Recently, the DR was derived by Maseso and Martellini \cite{maeso2020maximizing} from the stochastic portfolio theory
under the name of excess growth rate.
The concept of DR as well as its relationships to other diversification concepts have been extensively discussed in finance literature, see
\cite{erb2006strategic, greene2011sources, qian2012diversification, chambers2014limitations, cuthbertson2015diversification, meyer2017rebalancing, lesyk2021role} and the references therein.
In particular, the extensive empirical results reported in \cite{maeso2020maximizing} shows that maximizing DR of a portfolio may lead to
strong out-of-sample performance and hence results in competitive strategy 
in portfolio construction in certain market conditions. 
However, compared with the Markowitz mean-variance model
DR portfolios lack optimality theory for a deep understanding of their strong performance. This paper is to conduct a systematic theoretical study about DR portfolios and relate the optimal DR portfolios to the efficient frontier of
the mean-variance model. We will also illustrate the findings by DAX30 Index portfolios.

To motivate our research, let us consider the three fundamental results
that have become standard textbook material on modern
portfolio theory \cite{elton2009modern} about the mean-variance (MV) model. 
The first result is the {\bf efficient frontier} (EF) consisting of the portfolios that have the largest returns when the model is parametrized with the 
standard deviation. The EF starts from the Minimum Variance Portfolio (MVP) and
forms a smooth concave curve in the standard deviation and return space. 
The second result is the {\bf separation theorem}, which says that any portfolio on the EF can be represented by any two efficient portfolios. 
For example, MVP and the tangential portfolio are enough to produce the whole
EF. 
The last result is the {\bf capital market line} (CML) when there is a risk-free asset available to invest.  CML yields the largest Sharpe ratio. 
Moreover, to improve the performance of the mean-variance model, various constraints on the portfolio weights are often added to the model, see
\cite{jagannathan2003risk, demiguel2009generalized}.

Given DR is one kind of return, we may seek the highest DR of a portfolio when
it is parametrized along its standard deviation. In particular, we establish the
following results.

\begin{itemize}
	\item[(i)] The highest DR portfolios as a function of standard deviation form
	a concave and smooth curve in the standard deviation and DR space, which is
	denoted as $(\sigma, q)$ space with $q(\bfw)$ being the DR of a portfolio
	$\bfw$ (see Def.~\ref{Def-DR}). We refer to this curve as the efficient DR frontier.
	This concavity property is much like that of the efficient frontier of the mean-variance model in the standard deviation and expected return space.
	A distinctive feature is that the efficient DR curve is strongly concave.
	Consequently, it has the highest DR along the curve and it is the
	much studied Maximum Diversification Return Portfolio (MDRP). 
	
	\item[(ii)] There is also a separation theorem: any efficient DR portfolio
	is a simple convex combination of the two portfolios: MVP and MDRP.
	From the perspective of DR principle (the higher DR of a portfolio the better),
	any portfolio beyond MDRP is discarded as it would have lower DR, but higher risk.
	This is all due to the strong concavity of the efficient DR curve.
	
	\item[(iii)] When there is a risk-free asset available, the efficient DR frontier becomes a standard parabola and is also strongly concave.
	This means that a line similar to the CML in the $(\sigma, q)$ space does not exist.
	 
\end{itemize}

Furthermore, we investigate what the MV efficient portfolios would look like when putting in the $(\sigma, q)$ space.
We call the resulting curve the DR curve of MV efficient portfolios.
It turns out that under certain conditions, it is also strongly concave. 
In this case, a portfolio that has the highest DR on this curve exists and it is
called the Q-portfolio. It is interesting to note 
that the risk of the
Q-portfolio must be less than that of MDRP.
This property makes the Q-portfolio useful because otherwise MDRP would be preferred. 
We will see in the numerical part, Q-portfolio belongs to a cluster of portfolios, 
whose performance is close to that of a market portfolio.
In the special case when the variance vector of all assets are proportional to the vector of their expected returns, 
the efficient DR curve and the DR curve of the efficient portfolios become the same. Consequently, Q-portfolio becomes MDRP.

\begin{figure}[h!]
	\centering
	\includegraphics[width= 1\textwidth]{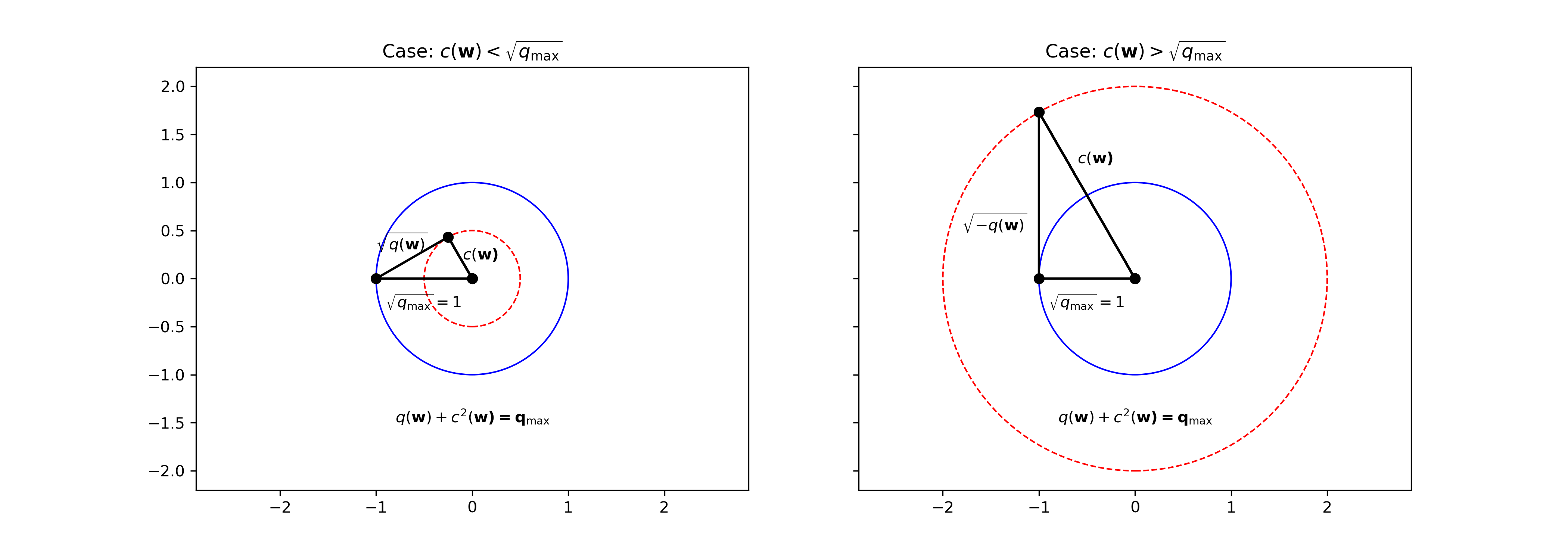}
	\caption{Pythagoras-style relationship between $c^2(\bfw)$ and $q(\bfw)$:
		the total sum is a constant. 
		The blue circle represents embedding sphere of the MDRP and the origin is MDRP 
		(see Thm.~\ref{Thm-Qi} for the detail). 
		Its radius is $\sqrt{q_{\max}}$ and it is one in this illustration.
		The red circle represents those portfolios that have the same centrality $c(\bfw)$.
		When the red circle is inside the blue circle, the length of the hypotenuse is $\sqrt{q_{\max}}$. 
		Otherwise, the length of the hypotenuse is $c(\bfw)$ in the right-angle triangle.
	}
	\label{Fig-Pythogoras}
	%\vspace{-16pt}
\end{figure}

%%%%%%%%%%%%%%%%%%%%%%%%%%%%%%%%%%%%%%%%%%%%%%
Some of the results, the separation theorem in particular, heavily depend on
the new geometric interpretation of DR portfolios recently studied in
\cite{qi2022geometric}, where the concept of portfolio centrality, denoted by
$c(\bfw)$, was introduced. In this paper, we establish an important identity
\be \label{qcw}
  q(\bfw) + c^2(\bfw) = q_{\max}, \ \ \ \mbox{for any portfolio} \ \bfw ,
\ee
where $q_{\max}$ is the DR attained by MDRP.
This Pythagoras-style relationship once again confirms that
the concept of DR is Euclidean.
If we use the analogy to conservation energy,
we may think the investment universe has a total and constant energy $q_{\max}$
while the potential energy $c^2(\bfw)$ and the kinetic energy 
$q(\bfw)$ can be transferred to each other.
In our case, $q(\bfw)$ may be negative.
This happens when $c^2(\bfw)$ is bigger than $q_{\max}$.
This possibility is illustrated in Fig.~\ref{Fig-Pythogoras}.

Moreover, the centrality is related to a weighted norm of $\bfw$. 
This link between the DR of a portfolio and its weighted norm 
brings out a new perspective of DR based portfolios.
To put it another way, norm-weighted portfolios 
\cite{demiguel2009generalized, gotoh2011role}
implicitly impose a condition on the diversification return. 
Hence, the much studied norm-regularized portfolios \cite{demiguel2009generalized, gotoh2011role} 
provide another strong motivation for a systematic 
study of optimal DR portfolio, echoing the empirical research on MDRP in \cite{maeso2020maximizing}. 

The maximum diversification ratio portfolio (MDP) of Choueufaty and Coignard \cite{choueifaty2008toward},
which is behind several FSTE 100 products \cite{carmichael2018rao},
is often confused with MDRP (their names are very similar).
In fact, they are closely related.
For the long-only case (portfolio weights are nonnegative),
the gap between the objectives of the two portfolios can be quantitatively bounded and the gap
can be small. This is well illustrated in the numerical part, where MDP is shown to be
one of clustered portfolios.

%%%%%%%%%%%%%%%%%%%%%%%%%%%%%%%%%%%%%%%%%%
The paper is organized as follows.
The next section formally introduces the DR of a portfolio and provides
the new characterization (\ref{qcw}) relating the DR to the centrality of the portfolio, see Prop.~\ref{Prop-Entrop-Invariance}.
In Section~\ref{Section-EF-DRP}, we study the most efficient DR portfolio
under a given risk level and characterize the efficient DR curve as 
a strongly concave function, see Prop.~\ref{Prop-DR-Curve}.
We further prove the separation theorem (Thm.~\ref{Thm-Two-EF}) that states that
any portfolio on the efficient DR curve can be represented by a simple mixture of MVP and MDRP. 
In Section~\ref{Section-Two-EFs}, we study the DR representation of 
Markowitz efficient portfolios and characterize its relationship with
the efficient DR curve. 
The relationship is graphically represented in Fig.~\ref{fig:DR-Frontier},
very similar to  the classical graphical representation of the efficient 
frontier for the mean-variance model.
We also extend this comparison study to the case of including a risk-free asset.
We demonstrate the behaviour of the DR portfolios for DAX30 Index stocks in
Section~\ref{Section-Numerical}.
The paper concludes in Section~\ref{Section-Conclusion}.\\

%%%%%%%%%%%%%%%%%%%%%%%%%%%%%%%%%%
{\bf Notation:} Through the paper, a boldfaced lower letter denotes a column vector. For example, $\bfw \in \Re^n$ is a column vector of dimension $n$.
Its transpose $\bfw^\top$ is a row vector. 
Let
$V$ denotes the covariance matrix of 
the returns of $n$ assets. 
Let $\eta$ denote its diagonal (column) vector of $V$, i.e., $\eta := \diag(V)$, where ``$:=$'' means ``define''.
From $V$, we define a new matrix:
\[
    D := \frac 12 \Big( \eta^\top \bfone_n + \bfone_n \eta^\top   \Big) - V,
\]
where $\bfone_n$ is the column vector of all ones of dimension $n$.
The $n \times n$ identity matrix is denoted by $I_n$.
We will see in the next section that the properties of $D$ play an important role in
our analysis. When $V$ is nonsingular, so is $D$ by \cite[Lemma~2.3]{qi2022geometric}. In this case, we let $\bfw_{\mvp}$ denote
the minimum variance portfolio
\[
   \bfw_{\mvp} := \frac{V^{-1} \bfone_n}{ \bfone_n^\top V^{-1} \bfone_n  }
   \quad
   \mbox{and} \quad
   \sigma^2_{\mvp} = \frac{1}{\bfone_n^\top V^{-1} \bfone_n },
  % a := \bfone_n^\top V^{-1} \bfone_n = \frac 1{\sigma^2_{\mvp}}, \quad
  % b := \bfone_n^\top V^{-1} \overline{\bfr},
\]
where $\sigma^2_{\mvp}$ is the variance of MVP.
% and $\overline{\bfr}$
%is the expected return vector of the $n$ assets. 
Throughout, $\| \cdot\|$ denotes the Euclidean norm.

%%%%%%%%%%%%%%%%%%%%%%%%%%%%%%%%%%%%%%%%%%%%%%%%%%%%%%%%%%%
\section{New Perspective of Diversification Return} \label{Section-DR}

Let us first recall the definition of DR of a portfolio.

%%%%%%%%%%%%%%%%%%%%%%%%%%%%%%%%%%%%%%%%%%%%%%%%%%%%
\begin{definition} \label{Def-DR}
	\cite{booth1992diversification}
	Suppose there are $n$ assets whose covariance matrix of their returns is denoted by $V$. The diversification return
	of a portfolio $\bfw$ satisfying the
	budget constraint $\bfone_n^\top \bfw =1$ is defined by
	\[
	q_V(\bfw) := \frac 12 \Big( \eta^\top \bfw - \bfw^\top V \bfw  \Big)
	\qquad \mbox{with} \ \ \eta := \diag(V).
	\]
	We often drop the dependence of $q$ on $V$ for simplicity.
\end{definition}

\noindent 
In words, the diversification return is half of the difference
between the weighted-average variance of the assets
in the portfolio and the portfolio variance, see also \cite[Eq.~(9)]{maeso2020maximizing}.
We note that the factor $1/2$ in the definition of the diversification return is
important. Mathematically, it results from the second-order Taylor expansion of
the expected return of a portfolio, see \cite{booth1992diversification}.

In this part, we explain a peculiar phenomenon about DR. Different
representation of a portfolio may have different DR depending 
how the portfolio is decomposed. 
We then derive an identity relating the centrality of a portfolio to its DR. 
This identity will explain why and how DR may be negative. 

%%%%%%%%%%%%%%%%%%%%%%%%%%%%%%%%%%%
\subsection{The matter of size} \label{Subsection-size}

Let us use an example to demonstrate the fact that DR of a portfolio depends on
how many assets are in it. 

\begin{example} \label{Example-DR}
Suppose we have two independent risky assets $S_1$ and $S_2$, whose covariance matrix is the identity matrix. Suppose there is also a riskfree asset $S_0$. We consider a
portfolio $\bfw_m := (1/2, 1/4, 1/4)^\top$ investing $50\%$ of the wealth to $S_0$ and $25\%$ each to the two risky assets. The covariance matrix of the three
assets is
\[
    V = \left[
    \begin{array}{ccc}
    	0 & 0 & 0 \\
    	0 & 1 & 0 \\
    	0 & 0 & 1
    \end{array} 
    \right] \quad \mbox{with} \quad \eta = \left[  
     \begin{array}{c}
     	0 \\ 1 \\ 1
     \end{array} 
    \right] .
\]
The DR of the portfolio is $q_V(\bfw_m) = \frac 12( 1/2 - 1/8) = \frac{3}{16}$.

Now consider the equal-weight portfolio of the two risky asset $\bfw_{\ew} = (1/2, 1/2)^\top$. Its variance is $\sigma^2_{\ew} = 1/2$.
Then the portfolio $\bfw_m$ can be regarded as investing
$50\%$ to $S_0$ and $50\%$ to $\bfw_{\ew}$. 
The new representation of $\bfw_m$ and its covariance matrix of
the assets $S_0$ and $\bfw_{\ew}$ are given below:
\[
\bfw_m = \left(  \begin{array}{c}
	 1/2 \\
	 \frac 12 \bfw_{\ew}
	\end{array} \right) ,
\qquad 
  \widetilde{V} = \left[
  \begin{array}{cc}
  	0 & 0 \\
  	0 & 1/2
  \end{array} 
  \right] \quad \mbox{and} \quad \widetilde{\eta} = \left[
  \begin{array}{c}
  	0 \\ 1/2
  \end{array} 
  \right].
\]
Then the DR of $\bfw_m$ can also be calculated as 
$q_{\widetilde{V}} (\bfw_m) = \frac 12 ( 1/4 - 1/8) = \frac 1{16}$.
\end{example}

On the surface, it seems that two different DRs of the same portfolio $\bfw_m$ in
Example~\ref{Example-DR} are due to the fact that two different covariance matrices are being used. 
On a deep level, it is because the DR of a portfolio measures the difference between the expected return of the portfolio and its second-order approximation through
the expected returns of its constituents. Hence, the number of its constitutes is vital in 
computing DR. For $q_{V}(\bfw_m)$, there are three assets while $q_{\widetilde{V}}(\bfw_m)$ only used two (one is $S_0$ and the other is $\bfw_{\ew}$). 
Therefore, it is only meaningful to discuss DR after the number of the assets
involved has been decided.

%%%%%%%%%%%%%%%%%%%%%%%%%%%%%%%%%%
\subsection{Centrality of portfolio}

This part introduces an important concept called the centrality of portfolio and
study its relationship with the DR of the portfolio. 
%Suppose we have $n$ risky assets with $V$ being the covariance matrix of the
%returns of those assets. 
For any portfolio $\bfw$ satisfying the budget constraint $\bfone_n^\top \bfw =1$, its DR can be represented as follows:
\[
q(\bfw) = \frac 12 \eta^\top \bfw - \frac 12 \bfw^\top V \bfw
= \frac 12 \bfw^\top \underbrace{\Big(   
	(\eta \bfone_n^\top + \bfone_n \eta^\top )/ 2 - V
	\Big)}_{=\; D} \bfw
\]
Recall that $\eta = \diag(V)$. This implies that the diagonal of $D$ are zeros and
all its elements $D_{ij} \ge 0$. In fact, $D$ is Euclidean Distance Matrix (EDM)
\cite[Lemma~2.2]{qi2022geometric}. 
This means that there exist a set of points $\bfx_i$, $i=1, \ldots, n$ in the
Euclidean space $\Re^k$ for some $k>0$ such that
\[
  \| \bfx_i - \bfx_j \|^2 = D_{ij}, \qquad i, j=1, \ldots, n.
\]
Those points are called the embedding points of $D$ because their pairwise
squared Euclidean distances recover the elements in $D$. 
This embedding result can be traced back to Schoenberg \cite{schoenberg1935remarks} and
Young and Householder \cite{young1938discussion}.
Obviously, there exist infinitely many such embedding points because
shift and rotation transformations do not change pairwise distances. 
There is one special set of embedding points that are useful in characterizing 
the MDRP. We explain how we derive those points. 

In \cite[Thm.~3.1]{qi2022geometric}, it is showed that
\[
  \bfw_{\mdrp} = \frac{D^- \bfone_n}{ \bfone_n^\top D^- \bfone_n },
\]
where $D^-$ can be chosen as any generalized inverse of $D$ satisfying $DD^-D =D$
and $D^-DD^-=D^-$. That is, $\bfw_{\mdrp}$ does not depend on the choice of $D^-$.
Let
\be \label{Js}
\bfs := \bfw_{\mdrp}, \quad 
J_{\bfs} := I_n - \bfs \bfone_n^\top \quad \mbox{and} \quad
B := -\frac 12 J_{\bfs}^\top D J_{\bfs}.
\ee
According to the theory of EDM \cite{gower1982euclidean}, the matrix $B$
is positive semidefinite. Suppose it has the following spectral
decomposition:
\be \label{B-Decomp}
   B = [\bfp_1, \ldots, \bfp_k] \left[
   \begin{array}{ccc}
   	 \lambda_1 &  & \\
   	  & \ddots & \\
   	  & & \lambda_k 
   	\end{array} 
   \right] \left[
   \begin{array}{c}
   	 \bfp_1^\top \\
   	 \vdots \\
   	 \bfp_k^\top
   	\end{array} 
   \right],
\ee
where $\lambda_1 \ge \cdots \ge \lambda_k >0$ are the positive eigenvalues of
$B$ and $\bfp_i \in \Re^n$ are the corresponding orthonormal eigenvectors. Let
\be \label{Xpoints}
  X := [ \bfx_1, \ldots, \bfx_n] = \left[
  \begin{array}{ccc}
  	\sqrt{\lambda_1} &  & \\
  	& \ddots & \\
  	& & \sqrt{\lambda_k }
  \end{array} 
  \right] \left[
  \begin{array}{c}
  	\bfp_1^\top \\
  	\vdots \\
  	\bfp_k^\top
  \end{array} 
  \right] .
\ee
Then $B=X^\top X$ and $\bfx_i$, $i=1, \ldots, n$ are one set of embedding points of $D$. Moreover, the following are proved in \cite{qi2022geometric}.

\begin{theorem} \label{Thm-Qi}
\cite[Thm.~3.1]{qi2022geometric}	
Let $q_{\max} := q(\bfw_{mdrp})$ be the maximum diversification return of MDRP.
Let the embedding points $\bfx_i$, $i=1, \ldots, n$ be obtained by (\ref{Xpoints}). 
Then $q_{\max} = \frac{1}{2 \bfone_n^\top D^- \bfone_n} >0$ and  those embedding points sit on the sphere of centered at the origin with the radius
$ \sqrt{q_{\max}}$, i.e.,
\[
  \| \bfx_i \| = \sqrt{ \frac{1}{2 \bfone_n^\top D^- \bfone_n}  }
               = \sqrt{q_{\max}}, \qquad i=1, \ldots, n.
\]
Furthermore, the origin is the $\bfs$-weighted center of the embedding points:
\[
   0 = s_1 \bfx_1 + s_2 \bfx_2 + \cdots + s_n \bfx_n  \quad \mbox{or equivalently} \quad X\bfs = 0.
\]
\end{theorem}

Geometrically, Thm.~\ref{Thm-Qi} means that the origin represents MDRP. 
For any other portfolio $\bfw$, the
$\bfw$-weighted centre of the embedding points is away from the origin.
Its distance from the origin measures how far it is from MDRP. 
We call the distance portfolio centrality.

\begin{definition} \label{Def-Centrality}
	(Portfolio centrality)
%Let $V$ be the covariance matrix of the returns of $n$ assets and its corresponding Euclidean distance matrix $D = 
% (\eta\bfone_n^\top + \bfone_n \eta^\top)/2 - V$, where $\eta \in \diag(V)$. 
Let the matrix $B$ be defined by (\ref{Js}).
For a given portfolio $\bfw$, its centrality is defined by
\[
  c(\bfw) := \sqrt{\bfw^\top B \bfw  }.
\]
%and the corresponding embedding points $\{ \bfx_i \}_{i=1}^n$ be constructed by (\ref{Xpoints}). 
%For a given portfolio $\bfw$, the quantity $c(\bfw)$ below is called the
%centrality of portfolio $\bfw$:
%\[
%   c(\bfw) = \| w_1 \bfx_1 + w_2 \bfx_2 + \cdots + w_n \bfx_n \| = \| X \bfw \| .
%\]
\end{definition}

\begin{remark} \label{Remark-Centrality}
(Interpretation of centrality)
Let us explain why $c(\bfw)$ measures the centrality of portfolio $\bfw$ with respect to
MDRP. Using (\ref{Xpoints}), we have
\[
  c^2(\bfw) = \bfw^\top B \bfw = \bfw^\top X^\top X \bfw
  = \| X \bfw \|^2 .
\]
Therefore, $c(\bfw) = \| X\bfw\|$. 
From Thm.~\ref{Thm-Qi}, we have $X\bfs =0$, which is the centre of the embedding sphere.
On the  other hand, $X\bfw$ represents the embedding point of portfolio $\bfw$ in the
embedding space spanned by $\{ \bfx_i\}$. 
The quantity $\| X \bfw\|$ 
measures how far it is from the origin. 
	
Another interesting interpretation is as follows. The decomposition
(\ref{B-Decomp}) gives rise to $k$ principle dimensions $\bfp_i$, $i=1, \ldots, k$. For a given portfolio $\bfw$, we may compute its principle coordinates
\[
   \widetilde{w}_i := \bfp_i^\top \bfw, \quad i=1, \ldots, n.
\]
Denote this new point by $\widetilde{\bfw}$. Obviously, the length remains the same: $\| \widetilde{\bfw}\| = \| \bfw\|$. 
The centrality considers the weighted length by the eigenvalues:
\[
  c(\bfw) = \sqrt{\lambda_1  \widetilde{w}_1^2 + \cdots + \lambda_k  \widetilde{w}_k^2 }.
\]
This interpretation is closely related to the principle coordinate analysis, 
initially studied by Gower \cite{gower1966some}. 

\end{remark}

Since the centrality measures how far a portfolio is from the origin and
the origin represents the maximum diversification return, 
we like to know if the higher centrality means lower diversification return. 
The following result proves it is the case.

\begin{proposition} \label{Prop-Entrop-Invariance}
	For any portfolio $\bfw$ satisfying $\bfone_n^\top \bfw = 1$, it holds
	\be \label{Centrality-Identity}
	c^2(\bfw) + q(\bfw) = q_{\max} .
	\ee
\end{proposition}

{\bf Proof.} 
We first note from Thm.~\ref{Thm-Qi} that
\[
q_{\max} = q(\bfs)
= \frac 12 \frac{1}{\bfone^n D^{-} \bfone_n  }.
\]
For any portfolio $\bfw$, it is easy to see
\[
J_{\bfs} \bfw = \Big( I_n - \bfs \bfone_n^\top  \Big) \bfw
= \bfw - \bfs   \qquad (\mbox{using } \bfone_n^\top \bfw = 1)
\]
and
\[
\bfw^\top D\bfs = \frac{ \bfw^\top D D^{-} \bfone_n}{\bfone_n^\top D^{-} \bfone } 
= \frac{ \bfw^\top \bfone_n}{\bfone_n^\top D^{-} \bfone } 
= \frac{1}{\bfone_n^\top D^{-} \bfone} = 2 q(\bfs),
\]
where the second equation used the fact $DD^-\bfone_n = \bfone_n$,
which was established in \cite[Thm.~2]{gower1985properties} for
any Euclidean distance matrix. 
Using those facts, we compute $c^2(\bfw)$ as follows:
\begin{eqnarray*}
	c^2(\bfw) &=& \bfw^\top B \bfw 
	= - \frac 12 \bfw^\top J_{\bfs}^\top D J_{\bfs} \bfw \\
	&=& - \frac 12 (\bfw- \bfs)^\top D (\bfw- \bfs)
	= -\frac 12 \bfw^\top D \bfw + \bfw^\top D \bfs - \frac 12 \bfs^\top D \bfs \\
	&=& - q(\bfw) + 2 q(\bfs) - q(\bfs) 
	= - q(\bfw) + q(\bfs).
\end{eqnarray*} 
Noticing $q(\bfs) = q_{\max}$, we derived the claimed identity. \hfill $\Box$ \\

The identity (\ref{Centrality-Identity}) has a significant implication:
for a portfolio, if its weighted center is further away from the origin, then
its centrality increases and consequently its diversification return decrease. 
The sum of the diversification return of a portfolio and its squared centrality
is a constant, which is the maximum diversification 
return $q_{\max}$. 
We recall that the diversification return can be cast as Rao's quadratic entropy \cite{rao1982diversity, qi2022geometric}, and entropy is generally one
kind of energy. In this sense, the considered $n$ risky assets have a constant
energy and the energy decomposes into two parts for any portfolio. 
One part is allocated to the diversification return of the portfolio and 
the other part is for its centrality (distance from the location of total energy). 
At any point of state, the two types of energy can transfer to each other.
But their total sum remains constant.
It becomes an art of portfolio construction to decide which state an investor would like to 
stay and invest.

%%%%%%%%%%%%%%%%%%%%%%%%%%%%%%%%%%%%%%%%%%%%%%%%%%%%%%
\subsection{Norm-regularized portfolios are diversified}

Since the seminar work of Jaganathan and Ma \cite{jagannathan2003risk},
norm constrained portfolio construction has gained significant attention, 
see, e.g., 
\cite{demiguel2009generalized, gotoh2011role}.
Our key message below is that norm constraints implicitly require
the diversification return to be above a certain level. 
Let us use the global minimum variance portfolio (GMVP) 
under a general norm constraint
studied in \cite{demiguel2009generalized} to demonstrate this implication.
\be \label{GMVP}
  \min \frac 12 \bfw^\top V \bfw , \quad \mbox{s.t.}
  \quad \bfone_n^\top \bfw =1, \ \ 
  \| \bfw \|_A \le \tau,
\ee 
for some $\tau>0$, where $\| \bfw\|_A$ is the A-norm of $\bfw$ defined by
$
 \| \bfw\|_A : = \sqrt{ \bfw^\top A \bfw }
$
with $A$ being symmetric and positive definite. In particular, when $A=I_n$ (the identity matrix), $\|\bfw\|_A = \| \bfw\|$ becomes the Euclidean norm (also known as
the $\ell_2$-norm in \cite{yen2014solving,zhao2021optimal}).

Recall the matrix $B$ is defined by (\ref{Js}).
Let 
$
  \beta^2 := \lambda_{\min}(A) /\lambda_{\max}(B) >0,
$
where $\lambda_{\max}(B)$ is the largest eigenvalue of $B$ and 
$\lambda_{\min}(A)$ is the smallest eigenvalue of $A$. 
Then the matrix $(A - \beta^2 B)$ is positive semidefinite. Consequently, we have
\[
  \| \bfw\|_A \ge \beta \| \bfw\|_B = \beta c(\bfw).
\]
Therefore, the norm constraint in (\ref{GMVP}) implies
\[
  c(\bfw) \le \tau/\beta,
\]
or equivalently by using the identity (\ref{Centrality-Identity})
\[
   q(\bfw) \ge q_{\max} - (\tau/\beta)^2 .
\]
This means that any $A$-norm constraint in (\ref{GMVP}) actually requires the
diversification return of the portfolio be above a threshold. 
Interestingly, Carmichael et. al. \cite{carmichael2018rao} imposed a diversification return constraint of the type $q(\bfw) \ge \tau$ in 
their portfolio construction and observed strong out-of-sample performance,
similar to what has been observed in \cite{demiguel2009generalized}. 
The connection between the norm constraint and the diversification return constraint explains why both can lead to similar out-of-sample performance.

%%%%%%%%%%%%%%%%%%%%%%%%%%%%%%%%%%%%%%%%%%%%%%%%%%%%%%%%%%
%Furthermore, if $V$ is assumed to be nonsingular, then we have the following
%equivalent characterization of $q(\bfw)$.
%
%\begin{proposition} \label{Prop-qA}
%Suppose $V$ is nonsingular. Let $A := B + \bfone_n \bfone_n^\top$, where $B$
%is defined in (\ref{Js}). Then $A$ is positive definite and
%\[
%   q(\bfw) = q_{\max} + 1 - \| \bfw\|_A^2 \qquad \mbox{for any portfolio} \
%   \bfw \ \mbox{satisfying} \ \bfone_n^\top \bfw = 1.
%\]
%\end{proposition}
%
%{\bf Proof.}
%It follows from \cite[Lemma~2.3]{qi2022geometric} that $D$ is nonsingular
%under the nonsingularity assumption of $V$. 
%Consequently, the rank of the matrix $B$ in (\ref{Js}) is $(n-1)$ because
%the rank of the matrix $J_{\bfs}$ is $(n-1)$.
%Matrix $B$ has only one zero eigenvalue. 
%In other words, the vector $\bfs$ is the eigenvector of $B$ corresponding to its
%zero eigenvalue. This implies the positive definiteness of the matrix $A$.
%We note that
%\[
% c^2(\bfw) = \bfw^\top B \bfw 
% = \bfw^\top\top A \bfw - (\bfone_n^\top \bfw)^2
% = \| \bfw\|_A^2 - 1.
%\]
%The result follows from the identity (\ref{Centrality-Identity}). 
%\hfill $\Box$ \\
%%%%%%%%%%%%%%%%%%%%%%%%%%%%%%%%%%%%%%%%%%%%%%%%%%%%%%%%

We now borrow an example from \cite{qi2022geometric} to illustrate the main
results in this section.

\begin{figure}[h!]
	\centering
	\includegraphics[width= 0.8\textwidth]{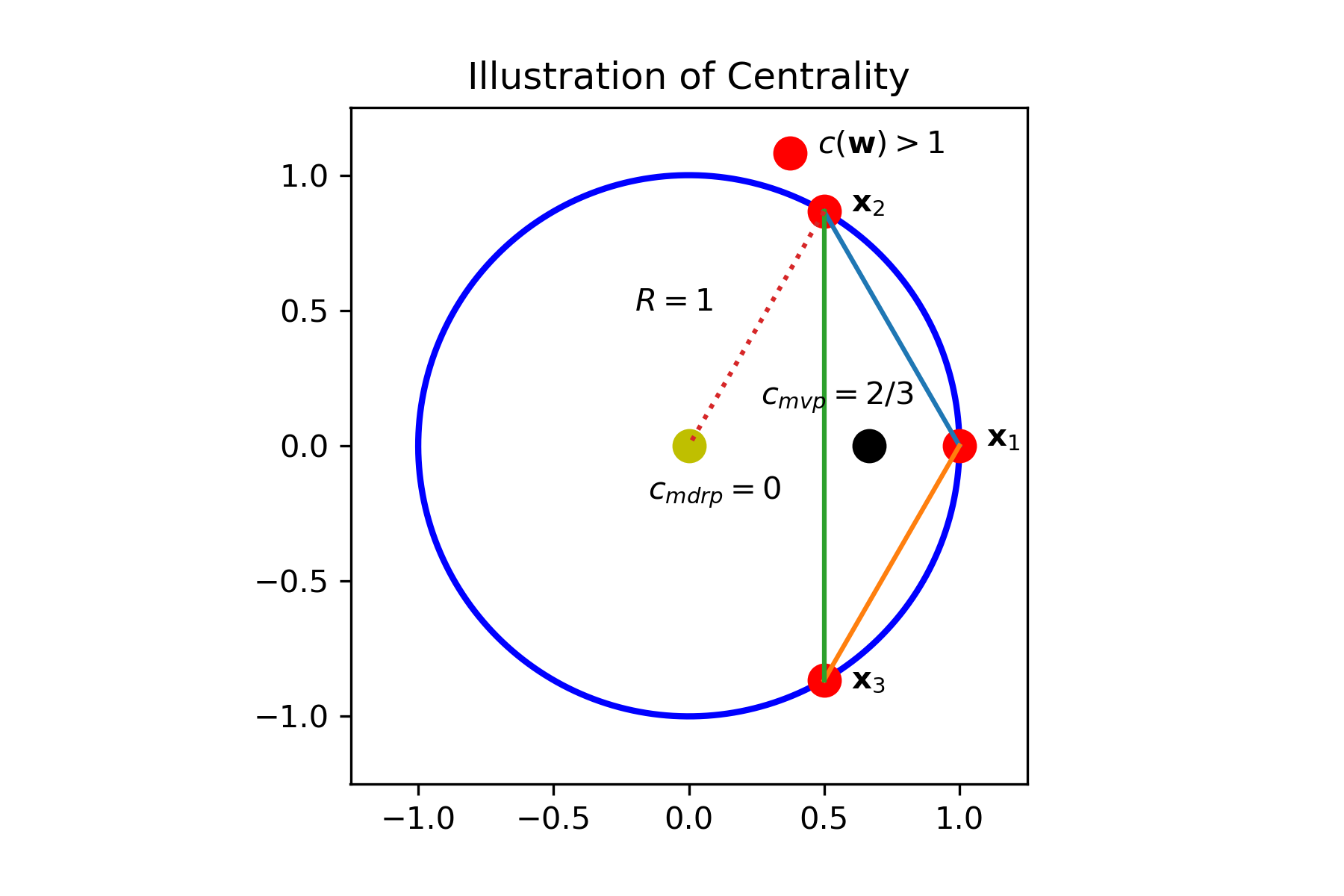}
	\caption{Asset and portfolio representation in the embedding space for the data in Example~\ref{Example-Centrality}: 
		the sphere is centred at the origin with $R=1$ and the three embedding 
		points on the sphere of the assets are $\bfx_1$, $\bfx_2$ and $\bfx_3$.
		Three portfolios are: MDRP (the origin), MVP (dot in black with its centrality being $2/3$), a portfolio outside of the sphere ($c(\bfw) >1$).
		}
	\label{Fig-Centrality}
	%\vspace{-16pt}
\end{figure}

\begin{example} \label{Example-Centrality}
	\cite[Example 3.1]{qi2022geometric}
Suppose there are three risky assets $S_i$, $i=1,2,3$, whose covariance matrix $V$, the corresponding
distance matrix $D$ and its inverse $D^{-1}$ are respectively given by
\[
V = \frac 19 \left[
\begin{array}{ccc}
	11 & 8 & 8 \\
	8 & 23 & -4 \\
	8 & -4 & 23
\end{array} 
\right], \qquad D =  \left[
\begin{array}{ccc}
	0 & 1 & 1 \\
	1 & 0 & 3 \\
	1 & 3 & 0
\end{array} 
\right], \qquad
D^{-1} =  \frac 16 \left[
\begin{array}{ccc}
	-9 & 3 & 3 \\
	3 & -1 & 1 \\
	3 & 1 & -1
\end{array} 
\right].
\]
We note that $V$ is also nonsingular. 
It follows from Thm.~\ref{Thm-Qi} that the MDRP is
$
\bfw_{\mdrp} = (-1, 1, 1)^\top.
$
The corresponding embedding points are
$
\bfx_1 = (1, 0)^\top,
$
$
\bfx_2 = (1/2, \sqrt{3}/2)^\top,
$
and
$
\bfx_3 = (1/2, -\sqrt{3}/2)^\top .
$
They lie on the sphere of radius $R =1$ centred at the origin, see
Fig.~\ref{Fig-Centrality}, where we also indicated the minimum variance 
portfolio $\bfw_{\mvp} = (1/3, 1/3, 1/3)^\top$
 with $c_{\mvp}=2/3$ or equivalently $q_{\mvp} = 1-(2/3)^2 = 5/9$.
Outside the unit circle, it holds $c(\bfw) >1$, which implies negative diversification 
returns for those portfolios. 
The triangle formed by $\bfx_1$, $\bfx_2$ and $\bfx_3$ is the region of portfolios with nonnegative weights  because the region is the convex hull of the three points. Outside of the triangle, there must be at least one negative weight. For example, MDRP is outside the triangle and its first weight is ($-1$). It is also straightforward to impose constraints on the diversification return $q(\bfw)$. For example, if we require the DR to be nonnegative, then
we have
\[
 q(\bfw) \ge 0  \quad \Longleftrightarrow \quad
 c^2(\bfw) \le 1 \quad \Longleftrightarrow \quad
 \| \bfw \|_B \le 1 \quad \Longleftrightarrow \quad
 \frac 14 (1+w_1)^2 + \frac 34 (w_2 - w_3)^2 \le 1.
\]

\end{example}

We will see the identity (\ref{Centrality-Identity}) greatly facilitate our derivation 
of the separation theorem below.

%%%%%%%%%%%%%%%%%%%%%%%%%%%%%%%%%%%%%%%%%%%%%%%%%%%%%%%%%%%
\section{Efficient Frontier of DRP} \label{Section-EF-DRP}

In the proceeding section, we established the identity (\ref{Centrality-Identity}), which means that maximizing the diversification 
return is equivalent to minimizing the centrality, which in turn is closely
related to norm-regularization of the portfolio weight. 
In this part, we study the maximum DR portfolio under a given risk level $\sigma^2$ and then establish a separation theorem.
From now on, we assume that $V$ is nonsingular. Consequently, the Euclidean distance matrix $D$ is also nonsingular. Moreover, the maximum diversification return portfolio has the following closed-form solution: 
\be \label{w_mdrp}
\bfw_{\mdrp} = \frac{D^{-1} \bfone_n }{\bfone_n^\top D^{-1} \bfone_n }
= \left(
1 - \frac{\bfone_n^\top V^{-1} \eta}{2}
\right) \frac{V^{-1} \bfone_n}{\bfone_n^\top V^{-1} \bfone_n } + \frac 12 V^{-1} \eta.
\ee
The first representation above is in terms of $D^{-1}$, and the second is in terms of $V^{-1}$ and is also obtained in \cite[Eq~(18)]{maeso2020maximizing}.

%%%%%%%%%%%%%%%%%%%%%%%%%%%%%%%
\subsection{The efficient DR curve}

The optimization problem we study is defined as follows:
\be \label{DRP-Risk}
  \bfw(\sigma) := \arg\max q(\bfw), \quad \mbox{s.t.} \quad
  \bfone_n^\top \bfw = 1, \ \ 
  \bfw^\top V \bfw = \sigma^2.
\ee 
In particular, we study the relationship between $q(\bfw(\sigma))$ and
$\sigma$ so that the optimal portfolios $\bfw(\sigma)$ 
form a dominating smooth and concave curve 
in the $(\sigma, q)$ diagram. 
This development is much in the spirit of Merton's classical contribution
\cite{merton1972analytic} to the mean-variance model. 

Under the condition $\bfone_n^\top \bfw =1$, the objective $q(\bfw)$ becomes
\be \label{qw}
  q(\bfw) = \frac 12 \bfw^\top D \bfw = \frac 12 \eta^\top \bfw - \frac 12 \bfw^\top V \bfw 
  = \frac 12 \eta^\top \bfw - \frac 12 \sigma^2. 
\ee
Problem (\ref{DRP-Risk}) reduces to
\be \label{DRP-Risk-2}
 \bfw(\sigma) := \arg\max \frac 12 \eta^\top \bfw, \quad \mbox{s.t.} \quad
 \bfone_n^\top \bfw = 1, \ \ 
 \bfw^\top V \bfw = \sigma^2.
\ee 
This is to maximize a linear function over the surface of an ellipsoid.
Therefore, the optimization can be done over the ellipsoid, leading to the
following convex optimization problem:
\be \label{DRP-Convex}
 \bfw(\sigma) := \arg\max 2 \eta^\top \bfw, \quad \mbox{s.t.} \quad
 \bfone_n^\top \bfw = 1, \ \ 
 \bfw^\top V \bfw \le  \sigma^2.
\ee 
We note the small change in the objective: we used $2 \eta^\top \bfw$ instead of
$(\eta^\top \bfw)/2$. This change is a positive scaling of the objective function and hence it won't change the optimal solution $\bfw(\sigma)$. 
The benefit of this change will become clear when we study the KKT condition of
(\ref{DRP-Convex}). 
We also note that it is necessary to require $\sigma^2 \ge \sigma_{\mvp}^2$ 
(the variance of the minimum variance portfolio (MVP)) for 
$\bfw(\sigma)$ to be well defined.
%We also note that when $\sigma^2 = \sigma_{\mvp}^2$, then the feasible region
%of (\ref{DRP-Convex}) has only one point, which is $\bfw_{\mvp}$.
%Hence, for this case, $\bfw(\sigma) = \bfw_{\mvp}$. 
%We only need to consider the case $\sigma > \sigma_{\mvp}$. 
We have the following result.

\begin{proposition} \label{Prop-DR-Curve}
	For $\sigma^2 \ge \sigma_{\mvp}^2$, we have
\be \label{eta-w}
 \eta^\top \bfw(\sigma) = 
%  \sqrt{\eta^\top V^{-1} \eta - (\bfone_n^\top V^{-1} \eta )^2/ \bfone_n^\top V^{-1} \bfone_n  } 
 \rho \sqrt{ \sigma^2 - \sigma_{\mvp}^2} 
  + \frac{ \bfone_n^\top V^{-1} \eta}{ \bfone_n^\top V^{-1} \bfone_n   },
  \ \ 
  \mbox{with} \ \
  \rho := \sqrt{\eta^\top V^{-1} \eta - (\bfone_n^\top V^{-1} \eta )^2/ \bfone_n^\top V^{-1} \bfone_n  } .
\ee
Moreover, the diversification return of $\bfw(\sigma)$, denoted by 
$q_{\efdr}(\sigma) := q(\bfw(\sigma))$ is given by
\be \label{q_efdr}
q_{\efdr} (\sigma) =
-\frac 12 \Big(  \sqrt{\sigma^2 - \sigma_{\mvp}^2  } - \rho/2  \Big)^2
+ \frac 18 \rho^2 
+ q_{\mvp} .
\ee 
In particular, when $\sigma = \sigma_{\mvp}$,  the DR of the minimum variance portfolio, denoted by $q_{\mvp}$ is given by
\be \label{q_mvp}
q_{\mvp} =  \frac 12 \Big(  \bfone_n^\top V^{-1} \eta - 1 \Big) \sigma_{\mvp}^2 .
\ee 

%In other words, $(\sigma, q_{\efdr})$ forms the efficient frontier of
%the diversification return (efdr) portfolios in the digram of standard deviation and diversification return.
\end{proposition}

{\bf Proof.}
For the case $\sigma^2 = \sigma_{\mvp}^2$, the feasible region of (\ref{DRP-Convex}) has only one point, which is $\bfw_{\mvp}$. 
The formula (\ref{eta-w}) is easily versified because $\bfw(\sigma) = \bfw_{\mvp}$.

We now consider the case $\sigma^2 > \sigma_{\mvp}^2$. For this case, 
we note that the Slater condition is satisfied.
One can verify that $\bfw_{\mvp}$ is in the relative interior of the feasible
region of (\ref{DRP-Convex}). Therefore, the KKT condition holds at $\bfw =\bfw(\sigma)$:
\be \label{KKT}
 -2\eta - 2\lambda \bfone_n + 2 \beta V \bfw = 0, \quad
 \bfone_n^\top \bfw = 1, \quad
 \beta \ge 0, \  \ \bfw^\top V \bfw \le  \sigma^2 \ \ \mbox{and} \ \ 
 \beta (\bfw^\top V \bfw -  \sigma^2) = 0.
\ee 
where $(2 \lambda)$ is the Lagrange multiplier corresponding to $\bfone_n^\top \bfw =1$ and $\beta \ge 0$ is the Lagrange multiplier for the ellipsoid constraint. We solve the KKT condition (\ref{KKT}) below.

We first remove the trivial case that $\eta$ is proportional to $\bfone_n$, i.e.,
$\eta = \gamma \bfone_n$ for some $\gamma>0$. In this case, the objective 
function is constant $2\gamma$ and hence any point in the feasible
region of (\ref{DRP-Convex}) is optimal. It is easy to verify that
the formula in (\ref{eta-w}) gives the correct value for $\eta^\top \bfw(\sigma)$. Therefore, we assume $\eta \not= \gamma \bfone_n$ for any $\gamma>0$. It follows from the first equation in (\ref{KKT}) that
\[
  \beta \bfw = V^{-1} ( \eta + \lambda \bfone_n) \not= 0
  \qquad (\mbox{because}\  \eta + \lambda \bfone_n \not=0, \ \forall \ \lambda \in \Re ).
\]
Hence $\beta \not=0$. By using the second equation in (\ref{KKT}), we get
\[
 \bfw = \frac 1\beta \Big(  V^{-1} \eta + \lambda V^{-1}\bfone_n  \Big)
 \quad \mbox{and} \quad
 \lambda = \frac{1}{ \bfone_n^\top V^{-1} \bfone_n} \beta 
 - \frac{\bfone_n^\top V^{-1} \eta}{ \bfone_n^\top V^{-1} \bfone_n} .
\]
Since $\bfw(\sigma)$ must be on the boundary of the ellipsoid, we substitute 
$\bfw$ into the equation $\bfw^\top V \bfw = \sigma^2$ to get after
simplification:
\[
  \left( \sigma^2 - \sigma_{\mvp}^2 \right) \beta^2
  = \eta^\top V^{-1} \eta - \frac{(\bfone_n^\top V^{-1} \eta)^2}{ \bfone_n^\top V^{-1} \bfone_n}
\]
By Cauchy-Schwartz inequality, we have 
$
(\bfone_n^\top V^{-1} \eta)^2 \le (\eta^\top V^{-1} \eta ) (\bfone_n^\top V^{-1} \bfone_n)
$. Therefore, the right-hand side of the above equation is nonnegative.
Given $\beta > 0$, we have
\[
  \beta = \sqrt{ \frac{\eta^\top V^{-1} \eta - (\bfone_n^\top V^{-1} \eta)^2 / \bfone_n^\top V^{-1} \bfone_n}{ \sigma^2 - \sigma_{\mvp}^2 }     }. 
\] 
We then can compute $\lambda$ in terms of $\beta$ and subsequently we get
\begin{eqnarray*}
	\eta^\top \bfw (\sigma) &=& \frac 1{\beta} 
	\Big( \eta^\top V^{-1} \eta + \lambda \eta^\top V^{-1} \bfone_n  \Big) \\
	&=& \frac 1{\beta} 
	\left( \eta^\top V^{-1} \eta + \frac{\bfone_n^\top V^{-1} \eta}{ \bfone_n^\top V^{-1} \bfone_n} \beta 
	- \frac{(\bfone_n^\top V^{-1} \eta)^2}{ \bfone_n^\top V^{-1} \bfone_n}  \right) \\
	&=& \frac 1{\beta} 
	\left( \eta^\top V^{-1} \eta - \frac{(\bfone_n^\top V^{-1} \eta)^2}{ \bfone_n^\top V^{-1} \bfone_n} \right) 
	 + \frac{\bfone_n^\top V^{-1} \eta}{ \bfone_n^\top V^{-1} \bfone_n} .
\end{eqnarray*}
Substituting $\beta$ into the last equation above yields the formula
(\ref{eta-w}). We proved it for all cases. 

It follows from \eqref{eta-w} that
\[
  \eta^\top \bfw_{\mvp} = \frac{ \bfone_n^\top V^{-1} \eta  }{ \bfone_n^\top V^{-1} \bfone_n }
  = (\bfone_n^\top V^{-1} \eta ) \sigma^2_{\mvp}
  \quad (\mbox{using the fact}  \
  \sigma^2_{\mvp} = \frac{ 1 }{ \bfone_n^\top V^{-1} \bfone_n }),
\]
which, together with \eqref{qw} yields
\[
  q_{\mvp} = -\frac 12 \sigma^2_{\mvp} + \frac 12 \eta^\top \bfw_{\mvp}
           = \frac 12 (\bfone_n^\top V^{-1} \eta -1 ) \sigma^2_{\mvp} .
\]
This proves \eqref{q_mvp}).
Equation (\ref{q_efdr}) is simple application of (\ref{eta-w}) by noticing
\begin{eqnarray*}
	q_{\efdr} (\sigma) &=& q(\bfw(\sigma))
	= -\frac 12 \sigma^2 + \frac 12 \eta^\top \bfw(\sigma) \\
	&=& - \frac 12 \sigma^2 + \frac 12 \rho \sqrt{\sigma^2 - \sigma_{\mvp}^2}
	+ \frac 12 \frac{\bfone_n^\top V^{-1} \eta}{ \bfone_n^\top V^{-1} \bfone_n} \\
	&=& -\frac 12 \Big(  \sqrt{\sigma^2 - \sigma_{\mvp}^2  } - \rho/2  \Big)^2
	+ \frac 18 \rho^2 
	+ \underbrace{\frac 12 (\bfone_n^\top V^{-1} \eta) \sigma^2_{\mvp} - \frac 12 \sigma^2_{\mvp}}_{=\; q_{\mvp} \ \mbox{by} \ \eqref{q_mvp} }.
\end{eqnarray*}
%This establishes (\ref{q_efdr}).
\hfill $\Box$\\

We make the following remark, which includes some useful identities.

\begin{remark} \label{Remark-QMDRP}
 It is easy to verify that the function
		\[
		  f(x) := \left(\sqrt{x^2 - x_0^2} - c \right)^2
		\] 
		is strongly convex for $x \ge x_0$ when $x_0 >0$ and $c \ge 0$.
		Consequently, the function $q_{\efdr}(\sigma)$ is strongly concave in $\sigma$ for $\sigma \ge \sigma_{\mvp}$. Solving the equation
		\be \label{sigma-mdrp}
		 0 = q'_{\efdr}(\sigma) 
		 = \frac{\rho}{2} \frac{\sigma}{\sqrt{\sigma^2 - \sigma_{\mvp}^2}  } - \sigma
		\ee
		yields the maximum DR. The corresponding portfolio is
		$\bfw_{\mdrp}$. The solution of (\ref{sigma-mdrp}) and
		the corresponding DR are given by
		\be \label{qmdrp}
		  \sigma_{\mdrp}^2 = \sigma_{\mvp}^2 + \frac{\rho^2}4
		  \quad \mbox{and} \quad
		  q_{\mdrp} = \frac 18 \rho^2 + q_{\mvp}.
		\ee 
\end{remark}

The concavity of $q_{\efdr}(\sigma)$ makes the efficient DR portfolios sit on
the frontier of the $(\sigma, q)$ diagram, much like the classical efficient frontier of Markowitz Mean-Variance portfolios.
Next, we are going to prove that the frontier can be obtained by combining the
two portfolios: MVP and MDRP.

%%%%%%%%%%%%%%%%%%%%%%%%%%%%%%%%%%%%%%%%%%%%%%%%%%%%%%%%%%%%
\subsection{Separation Theorem} \label{Subsection-Separation-Thm}

The main purpose of this part is to prove that in order to invest on
the efficient frontier of DR, one only needs to invest in two 
portfolios: Minimum Variance Portfolio and Maximum Diversification Return Portfolio via their simple combination by specifying the level $\alpha$,
\be \label{walpha}
  \bfw_\alpha
  = \alpha \bfw_{\mdrp} + (1-\alpha) \bfw_{\mvp} 
  = \bfw_{\mvp} + \alpha (\bfw_{\mdrp} - \bfw_{\mvp} ), \quad \alpha \ge 0.
\ee
In the classical theory of Markowitz efficient frontier, this is known as
the separation theorem \cite{tobin1958liquidity} and \cite[Thm.~5.2.1]{brugiere2020quantitative}. 
We establish the result by proving that the standard deviation ($\sigma$) and
the diversification return ($q$) curve formed by the combined portfolio $\bfw(\alpha)$ is the same as $q_{\efdr}$. 

\begin{theorem}
	Let the constant $\rho$ be defined in (\ref{eta-w}).
	The optimal diversification return portfolio $\bfw(\sigma)$ of Problem (\ref{DRP-Risk}) 
	is given by
	\[
	   \bfw(\sigma) = \alpha \bfw_{\mdrp} + (1- \alpha) \bfw_{\mvp}, \qquad 
	   \mbox{with} \quad \alpha := \frac 2{\rho} \sqrt{\sigma^2 - \sigma^2_{\mvp}  } .
	\] 
\end{theorem}

{\bf Proof.}
Let us consider the affine combination $\bfw_{\alpha}$ of the portfolios
$\bfw_{\mdrp}$ and $\bfw_{\mvp}$ in (\ref{walpha}). We will show once the variance level of $\bfw_{\alpha}$ is specified as $\sigma^2$, then $\alpha$ 
is uniquely decided and equals $(2/\rho) \sqrt{\sigma^2 - \sigma^2_{\mvp} }$ as stated in the result.

Let 
\be \label{d_vector}
   \bfd := \bfw_{\mdrp} - \bfw_{\mvp} = \frac 12 V^{-1} \eta 
   - \frac{\bfone_n^\top V^{-1} \eta }{2} \times 
     \frac{ V^{-1} \bfone_n }{\bfone_n^\top V^{-1} \bfone_n}
\ee
and $\sigma^2_{\bfd}$ be the variance of $\bfd$. It is easy to verify that
\begin{eqnarray}
  \sigma^2_{\bfd}  = \bfd^\top V \bfd
  &=& \left(   
   \frac 12 V^{-1} \eta 
  - \frac{\bfone_n^\top V^{-1} \eta }{2} \times 
  \frac{ V^{-1} \bfone_n }{\bfone_n^\top V^{-1} \bfone_n}
  \right)^\top \left(
   \frac 12  \eta 
  - \frac{\bfone_n^\top V^{-1} \eta }{2} \times 
  \frac{ \bfone_n }{\bfone_n^\top V^{-1} \bfone_n}
  \right) \nonumber \\
  &=& \frac 14 \rho^2  \label{sigma_d}
\end{eqnarray}
and 
\[
 \bfd^\top V \bfw_{\mvp} 
 = \left(   
 \frac 12 V^{-1} \eta 
 - \frac{\bfone_n^\top V^{-1} \eta }{2} \times 
 \frac{ V^{-1} \bfone_n }{\bfone_n^\top V^{-1} \bfone_n}
 \right)^\top \frac{\bfone_n}{\bfone_n^\top V^{-1} \bfone_n   }
 = 0.
\]
 Consequently, we have
\[
 \sigma^2 = \bfw_{\alpha}^\top V \bfw_{\alpha} 
 = \sigma_{\mvp}^2 + \alpha^2 \sigma^2_{\bfd}
 = \sigma_{\mvp}^2 + \frac 14 \rho^2 \alpha^2.
\]
Given $\alpha \ge 0$, we have
\be \label{alpha-eq}
 \alpha = \frac 2{\rho} \sqrt{\sigma^2 - \sigma^2_{\mvp}   } .
\ee
Having established (\ref{alpha-eq}), we now compute the DR  of
$\bfw_{\alpha}$ and its centrality.
For simplicity, let $q_{\alpha} := q(\bfw_{\alpha})$ and
$c_{\alpha} := c(\bfw_{\alpha})$. We compute $q_{\alpha}$ through
$c_{\alpha}$. It is important to note that the matrix $J_{\bfs}$ defined in
(\ref{Js}) has the following property
\[
 J_{\bfs} \bfw_{\mdrp} = J_{\bfs} \bfs = 0.
\]
This implies $B \bfw_{\mdrp} = 0$.
We then have
\begin{eqnarray}
  c^2_{\alpha} &=& \bfw_{\alpha}^\top B \bfw_{\alpha} \nonumber \\ [1ex]
  &=& \alpha^2 \bfw_{\mdrp}^\top B \bfw_{\mdrp} 
  + 2 \alpha (1-\alpha) \bfw_{\mvp}^\top B \bfw_{\mdrp}
  + (1-\alpha)^2 \bfw_{\mvp}^\top B \bfw_{\mvp}   \nonumber \\ [1ex]
  &=& (1-\alpha)^2 \bfw_{\mvp}^\top B \bfw_{\mvp}
   = (1-\alpha)^2 c^2_{\mvp} \label{cw2}
\end{eqnarray}
We now make use of the identity (\ref{Centrality-Identity}) twice. 
The first time is on $c_{\mvp}$: the identity implies 
\[
  c^2_{\mvp} + q_{\mvp} = q_{\mdrp}  \qquad \mbox{(note $q_{\mdrp}$ is $q_{\max}$)}.
\]
It follows from (\ref{qmdrp}) that
\be \label{cmvp2}
  c^2_{\mvp} = \frac 18 \rho^2 .
\ee
The second time to use the identity (\ref{Centrality-Identity}) is on
$c_{\alpha}$ to get
\[
  q_{\alpha} = q_{\mdrp} - c^2_{\alpha}
  \stackrel{(\ref{cw2})}{=} q_{\mdrp} - (1-\alpha)^2 c^2_{\mvp} 
  \stackrel{(\ref{cmvp2})}{=} q_{\mdrp} -\frac{\rho^2}8  (1-\alpha)^2 .
\]
Substituting $\alpha$ in (\ref{alpha-eq}) and $q_{\mdrp}$ in (\ref{qmdrp}) into
above equation yields (after some simplification)
\[
 q_\alpha = -\frac 12 \left( \sqrt{\sigma^2 - \sigma^2_{\mvp}} - \frac{\rho}2   \right)^2 + \frac 18 \rho^2
 + q_{\mvp}
\]
Using formula (\ref{q_mvp}) for $q_{\mvp}$, $q_{\alpha}$ is just
$q_{\efdr}$ in (\ref{q_efdr}). This is to say that
the $(\sigma, q)$ curve formed by the portfolio $\bfw_{\alpha}$ is
same as that of the optimal portfolio $\bfw(\sigma)$. The correspondences between $\alpha$ and $\sigma$ is given by (\ref{alpha-eq}).
This proves the theorem. \hfill $\Box$ \\

It is worth noting that the meaningful range for $\alpha$ is $0 \le \alpha \le 1$.
This is because when $\alpha >1$, the corresponding portfolio has larger standard  deviation due to
\eqref{alpha-eq} than MDRP (i.e., $\sigma > \sigma_{\mdrp})$, but with less DR.

%%%%%%%%%%%%%%%%%%%%%%%%%%%%%%%%%%%%%%%%%%%%%%%%%%%%%%%%%%%%
\section{Relationship of Two Efficient Frontiers} \label{Section-Two-EFs}

In this part, we answer the interesting question where the classical 
efficient frontier would fit in the $(\sigma, q)$ diagram
and how far it would be from the DR frontier just studied. 
We also investigate the portfolio that has the highest diversification return on the efficient frontier. 

%%%%%%%%%%%%%%%%%%%%%%%%%%%%%%%%%%%%
\subsection{The $(\sigma, q)$ curve of MV efficient portfolios}

Let us briefly review some key facts about the efficient portfolios in 
Markowitz's Mean-Variance model. We denote the expected return vector of the $n$
risky assets by $\overline{\bfr}$. For a given level of risk $\sigma^2$, the MV 
model seeks the optimal portfolio by solving the following problem:
%\[
%   \bfw_{\mu} := \arg\min \bfw^\top V \bfw, \qquad
%   \mbox{s.t.} \quad \bfone_n^\top \bfw = 1, \quad
%                     \overline{\bfr}^\top \bfw = \mu.
%\]
\be \label{MV-sigma}
\bfw_{\mv} (\sigma) := \arg\max\; \overline{\bfr}^\top \bfw , \qquad
\mbox{s.t.} \quad \bfone_n^\top \bfw = 1, \quad
\bfw^\top V \bfw \le \sigma^2.
\ee
To remove the trivial case, we assume that $\overline{\bfr}$ is not proportional
to $\bfone_n$. Otherwise, all risky assets have the same expected returns and
as a consequence that all portfolios have the same return.
The classical portfolio theory says that all the optimal portfolios 
$\bfw_{\mv}$ as $\sigma$ varies form a efficient frontier in the
standard deviation and return diagram. 
Moreover, those portfolios have the following representation \cite[Corollary.~4.2.1]{brugiere2020quantitative}:
\[
  \bfw_{\mv}(\sigma) = \bfw_{\mvp} + \alpha \bfw_o, \quad \mbox{with} \ \
  %\alpha = \frac{\mu - r_{\mvp}}{m_o},
  \alpha = \sqrt{\sigma^2 - \sigma_{mvp}^2 },
\]
where 
\be \label{ab}
  \bfw_o := \frac{V^{-1} \Big( \overline{\bfr} - \frac ba \bfone_n \Big)  }
  {\sqrt{ \Big( \overline{\bfr} - \frac ba \bfone_n \Big)^\top V^{-1} \Big( \overline{\bfr} - \frac ba \bfone_n \Big)   }   }, 
  \quad 
  a := \bfone_n^\top V^{-1} \bfone_n, 
  \quad
  b := \bfone_n^\top V^{-1} \overline{\bfr}.
\ee
In \cite{brugiere2020quantitative}, $\bfw_o$ is called a self-financing 
portfolio and is normalized because of the following facts:
\be \label{Wo}
   \bfone_n^\top \bfw_o = 0 \qquad \mbox{and} \qquad
   \bfw_o^\top V \bfw_o = 1.
\ee 
We denote the diversification return of $\bfw_{\mv}(\sigma)$ by $q_{\ef}(\sigma)$ (here $q_{\ef}$ means it is the diversification return on the efficient frontier of the mean-variance model).

We now calculate $q_{\ef}$ by using the above facts and the formula 
\eqref{qw}:
\begin{eqnarray}
	q_{\ef}(\sigma) &=& -\frac 12 \sigma^2 + \frac 12 \eta^\top \bfw_{\mv}(\sigma) \nonumber \\ 
&=& -\frac 12 \sigma^2 + \frac 12 \eta^\top \bfw_{\mvp}
 + \frac {\alpha}2 \eta^\top \bfw_o \nonumber \\
&=& -\frac 12 \sigma^2 + \frac 12 \sqrt{ \sigma^2 - \sigma^2_{\mvp} }
(\eta^\top \bfw_o) + \frac 12 \eta^\top \bfw_{\mvp} \nonumber \\
&=& -\frac 12 \Big(  \sqrt{\sigma^2 - \sigma_{\mvp}^2  } - \frac 12 \eta^\top \bfw_o  \Big)^2
+ \frac 18 (\eta^\top \bfw_o)^2 
- \frac 12 \sigma^2_{\mvp}
+ \frac 12 \eta^\top \bfw_{\mvp}  \label{Eq-qef}
\end{eqnarray}
Note that 
\[
\eta^\top \bfw_{\mvp} = \eta^\top \frac{V^{-1} \bfone_n}{\bfone_n^\top V^{-1} \bfone_n }
=
(\eta^\top V^{-1} \bfone_n) \sigma^2_{\mvp}
= 2 q_{\mvp} + \sigma^2_{\mvp},
\]
where the last equality above used the identity in \eqref{q_mvp}. 
Substituting it into \eqref{Eq-qef}) gives 
\be \label{qef}
q_{\ef}(\sigma) = -\frac 12 \Big(  \sqrt{\sigma^2 - \sigma_{\mvp}^2  } - \frac 12 \eta^\top \bfw_o  \Big)^2
+ \frac 18 (\eta^\top \bfw_o)^2 
+ q_{\mvp}.
\ee 
The function $q_{\ef}(\sigma)$ defines the $(\sigma, q)$ curve of the
efficient frontier. Its form is similar to that of $q_{\efdr}(\sigma)$. 
Since $q_{\efdr}$ is the dominating frontier among all portfolios, we must have
\[
   q_{\efdr}(\sigma) \ge q_{\ef}(\sigma) \qquad \forall \sigma \ge \sigma_{\mvp}.
\]
However, unlike $q_{\efdr}$ being concave, the shape of $q_{\ef}$ depends on the sign of $\eta^\top \bfw_o$. 
We also like to understand how far $q_{\ef}$ is from $q_{\efdr}$. 
We make it precise below.

%%%%%%%%%%%%%%%%%%%%%%%%%%%%%%%%%
\subsection{The $Q$-portfolio}

The quantity difference between $q_{\efdr}$ and $q_{\ef}$ is closely related 
to the $Q$ portfolio:
\[
  \bfw_Q := \bfw_{\mvp} + ( \eta^\top \bfw_o/2) \bfw_o.
\]
Since $\bfone^\top \bfw_o = 0$ by \eqref{Wo}, $\bfw_Q$ is a portfolio satisfying
the budget constraint. Let us consider the self-financing portfolio
\[
   \bfd_Q := \bfw_Q - \bfw_{\mdrp} =  \frac 12 (\eta^\top \bfw_o) \bfw_o - \bfd
   \quad \mbox{with} \quad \bfd =  \bfw_{\mdrp} - \bfw_{\mvp}.
\]
It follows from (\ref{d_vector}) and the fact $\bfone_n^\top \bfw_o = 0$ that
\[
\bfw_o^\top V \bfd = \frac 12 \eta^\top \bfw_o.
\]
Using \eqref{sigma_d}, we can compute the variance of $\bfd_Q$ by
\[
 0 \le  \sigma^2_{\bfd_Q} = \frac 14 \rho^2 - (\eta^\top \bfw_o) \bfw_o^\top V \bfd 
  + \frac 14 (\eta^\top \bfw_o)^2 \underbrace{ \bfw_o^\top V \bfw_o }_{=1}
  = \frac 14 \rho^2 - \frac 14 (\eta^\top \bfw_o)^2 .
\]
Therefore, we have
\be \label{eta-rho}
  | \eta^\top \bfw_o | \le \rho.
\ee

\begin{theorem} \label{Thm-Two-EF}
	The following results concerning the two efficient frontiers hold.
\begin{itemize}
	\item[(i)] The $(\sigma, q)$ curve of $q_{\efdr}$ dominates 
	that of $q_{\ef}$. Moreover, 
	\[
	  q_{\efdr} (\sigma) - q_{\ef}(\sigma) = \frac 12 ( \rho - \eta^\top \bfw_o) \sqrt{\sigma^2 - \sigma^2_{\mvp}   } \qquad \ \forall \
	  \sigma \ge \sigma_{\mvp}.
	\] 
	
	\item[(ii)] Suppose $\eta^\top \bfw_o \ge 0$. 
	Then $q_{\ef}(\sigma)$ is strongly concave and 	
	the Q-portfolio is the 
	efficient portfolio that has the highest diversification return.
	Moreover, we have
	$
	  \sigma_{\mdrp} \ge \sigma_Q,
	$
	where $\sigma_Q$ is the standard deviation of the Q-portfolio.
	
	\item[(iii)] Suppose $\eta^\top \bfw_o < 0$. Let
	\[
	  \tau_o := \sigma_{\mvp} \sqrt{ 1 
	  	+ \frac{(\eta^\top \bfw_o)^{4/3} }{ \sigma_{\mvp}^{2/3} }  }.
	\]
	Then $q_{\ef}(\sigma)$ is strictly decreasing.
	Moreover, $q_{\ef}(\sigma)$ is convex over the interval $[\sigma_{\mvp}, \tau_o]$ and concave over $[\tau_o, \infty)$.
	
	\item[(iv)] When the risk vector is proportional to the expected return vector, i.e., $\eta = \gamma \overline{\bfr}$ for some $\gamma>0$, we have
	$\rho = \eta^\top \bfw_o$. This implies $q_{\efdr}(\sigma) = q_{\ef}(\sigma)$ for all $\sigma$.
	
\end{itemize}
	
\end{theorem}

{\bf Proof.}
(i) follows from the direct subtraction of $q_{\ef}$ of \eqref{qef} from
$q_{\efdr}$ of \eqref{q_efdr}. The dominance is because $(\rho - \eta^\top \bfw_o) \ge 0$ due to \eqref{eta-rho}.
 
For (ii), we notice that the function
\[
  f(x) = -\frac 12 \Big( \sqrt{x^2-x_0^2} - c  \Big)^2
\] 
is strongly concave when $x \ge x_0$,  $x_0 \ge 0 $, and $c\ge 0$. 
Simple application with $x= \sigma$, $x_0 = \sigma_{\mvp}$, and $c = \eta^\top \bfw_o \ge 0$ implies that $q_{\ef}(\sigma)$ is strongly concave.
In this case, $q_{\ef}(\sigma)$ reaches its maximum when
\be \label{Q-Condition}
  \sqrt{\sigma^2 - \sigma^2_{\mvp}} = \frac 12 \eta^\top \bfw_0.
\ee
This corresponds to the Q-portfolio. Hence 
(replacing $\sigma$ by $\sigma_Q$ in \eqref{Q-Condition})
\[
  \sigma_Q^2 = \sigma^2_{\mvp} + \frac 14 (\eta^\top \bfw_o)
  \stackrel{\eqref{eta-rho}}{\le} \sigma^2_{\mvp} + \frac 14 \rho^2 
  \stackrel{\eqref{sigma-mdrp}}{=} \sigma^2_{\mdrp}.
\] 
This proves $\sigma_Q \le \sigma_{\mdrp}$.

For (iii), we simply differentiate the function $q_{\ef}(\sigma)$ twice.
It is easy to see that its first derivative is always negative. Hence,
$q_{\ef}(\sigma)$ is decreasing. The second derivative is non-negative over the
interval $[\sigma_{\mvp}, \tau_o]$ and non-positive over $[\tau_o, \infty)$.

The claim in (iv) can be directly verified.
\hfill $\Box$ \\

%%%%%%%%%%%%%%%%%%%%%%%%%%%%%%%%%%%%%%%%%%%%%%%%%%%%
%\subsection{Graphical representation}

\begin{figure}
	\centering
	\begin{subfigure}{.5\textwidth}
		\centering
		\includegraphics[width=1\linewidth]{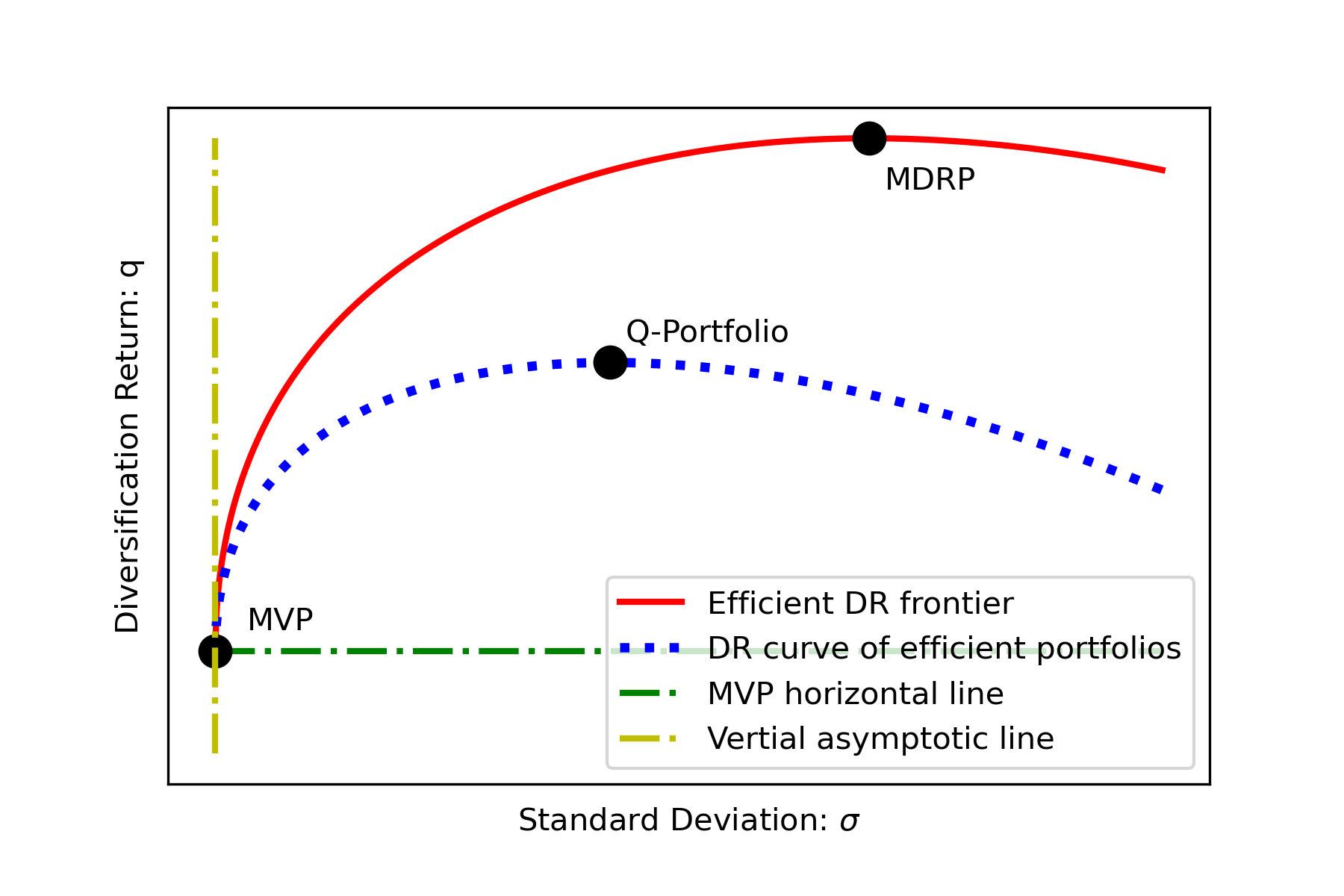}
		%  \caption{MVU}
		% \caption{Standard deviation vs diversification return.}
		\caption{}
		\label{fig:PositiveCase}
	\end{subfigure}%
	\begin{subfigure}{.5\textwidth}
		\centering
		\includegraphics[width=1\linewidth]{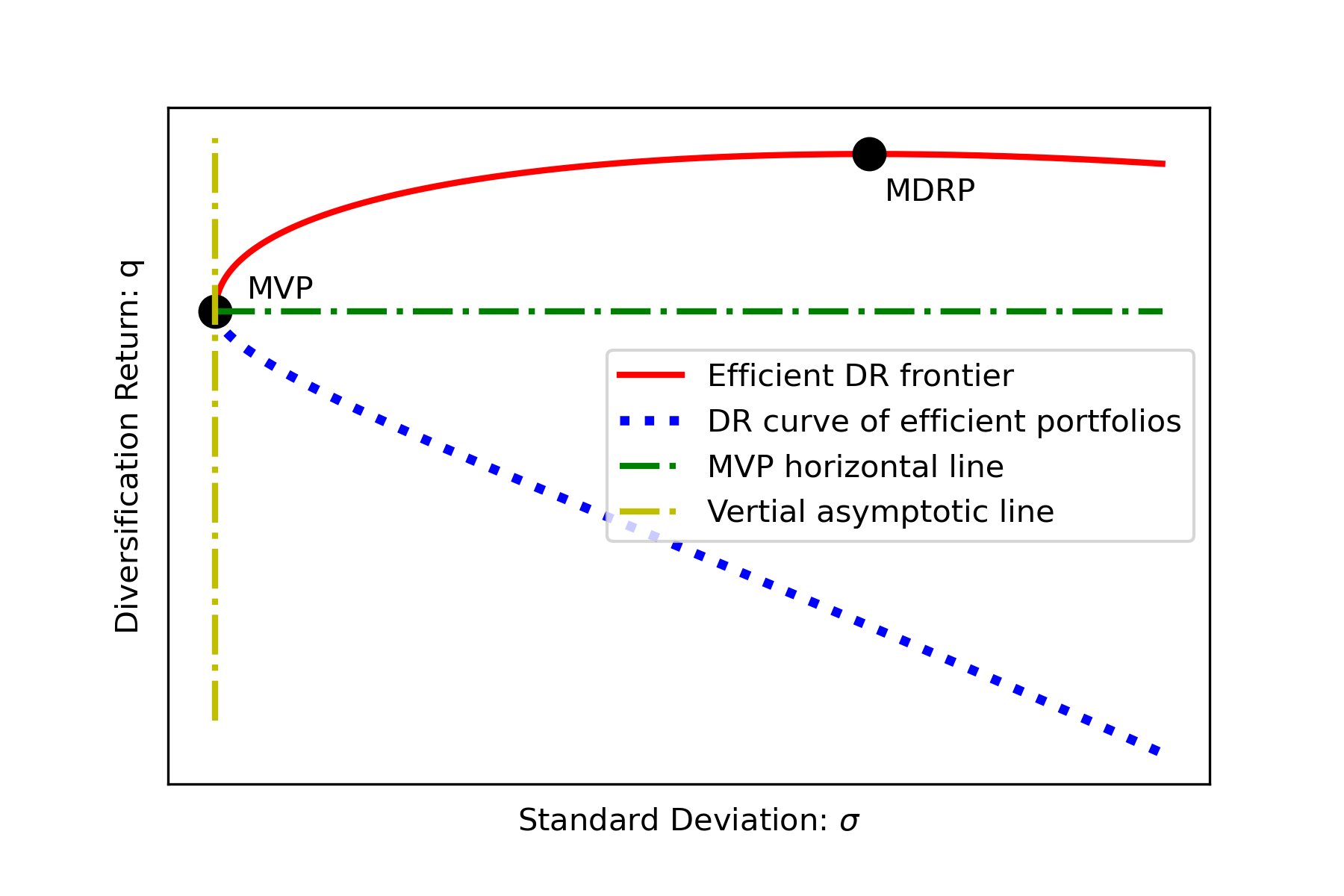}
		%\caption{Standard deviation vs expected return}
		\caption{}
		\label{fig:NegativeCase}
	\end{subfigure}
	\caption{Graphical representation of the efficient DR frontier in the $(\sigma, q)$ space.
		Fig.~\ref{fig:PositiveCase} is for the case when $\eta^\top \bfw_o \ge 0$.
		The Q-portfolio is on the curve of the efficient portfolios and has the 
		highest diversification return. It is below and on the left of MDRP.
		Fig.~\ref{fig:NegativeCase} is for $\eta^\top \bfw_o < 0$.
		The efficient portfolio that has the highest diversification return is MVP.
		In both figures, the efficient DR frontier dominates the DR curve of MV efficient portfolios.}
	\label{fig:DR-Frontier}
\end{figure}

We make the following remark regarding the relationship between the two DR curves.
The behaviour of the two DR curves and the key portfolios (MVP, MDRP, Q-portfolio) are illustrated in Fig.~\ref{fig:DR-Frontier} for both cases
($\eta^\top \bfw_o\ge 0$ or $<0$).

\begin{remark} \label{Remark-Two-EF}
	\begin{itemize}
	\item[(i)] The result in Thm.~\ref{Thm-Two-EF}(i) says that the risk-adjusted difference between the two curves is constant:
	\[
	  \frac{ q_{\efdr} (\sigma) - q_{\ef}(\sigma) }{ \sqrt{\sigma^2 - \sigma^2_{\mvp}   }}
	  = \frac 12 (\rho - \eta^\top\bfw_o), \quad \forall \ \sigma > \sigma_{\mvp}.
	\]
	This suggests that in order to retain the best of two worlds (MV efficiency and DR efficiency), one
	should focus on the region near the minimum variance portfolio. Too far away from it, the two efficient curves
	may go opposite directions and the gap grows at a constant rate in relation to the adjusted risk.
	In particular, when $\sigma_{\mvp}$ is small and hence the denominator is approximately proportional to $\sigma$, the gap is proportional to $\sigma$. 
	%That is an approximately linear difference.
	
	\item[(iii)] Under the condition $\eta^\top \bfw_o \ge 0$, the Q-portfolio 
	has the property that it is efficient and it has the highest DR among all
	MV efficient portfolios. Moreover, its standard deviation must be strictly less
	than the standard deviation of MDRP unless the Q-portfolio happens to be 
	MDRP. This is because $\bfd_Q \not=0$ if $Q$ portfolio is not MDRP and 
	$\sigma_{\bfd_Q} >0$, forcing $|\eta^\top \bfw_o| $ to be strictly less than $\rho$ in \eqref{eta-rho}. This in turn implies $\sigma_Q < \sigma_{\mdrp}$. If $\eta^\top \bfw_o <0$, the MVP is the efficient portfolio that has the highest DR among all MV efficient portfolios. 
	
	\item[(iv)] The result in Thm.~\ref{Thm-Two-EF}(iv) generalizes the well known fact that when the variance vector is proportional to the expected return vector, the maximum diversification return portfolio $\bfw_{\mdrp}$ is efficient,
	see \cite{carmichael2018rao, maeso2020maximizing}. The result says that any portfolio on the DR frontier is efficient. 
	
	\end{itemize} 
\end{remark}

%%%%%%%%%%%%%%%%%%%%%%%%%%%%%%%%%%%%%%%%%%%%%%%%%%%%%%%%%%

%%%%%%%%%%%%%%%%%%%%%%%%%%%%%%%%%%%%%%%%%%%%%%%%%%%%%%%%%%%
\subsection{Comparison when there is a risk-free asset}

Our study above can be straightforwardly extended to the important case when there is a risk-free asset. The classical Markowitz theory says the efficient portfolios form
the Capital Market Line (CML) \cite{sharpe1964capital, fama2004capital}.
We investigate how CML would look like on the $(\sigma, q)$ diagram. It turns out that
it is a parabola with some nice features.

Suppose the risk-free asset has a return $r_0$ and it is treated as the $(n+1)$th assets, appended to the $n$ risky assets studied above.
Let $\bfw_T$ denote the tangential portfolio (also known as the market portfolio). Then it is known (see e.g., \cite[Section 5.2]{brugiere2020quantitative}):
\[
  \bfw_T = \frac{1}{b- r_0 a} V^{-1} \Big(  \overline{\bfr} - r_0 \bfone_n  \Big) ,
\]
where the quantities $a$ and $b$ are defined in \eqref{ab}.
For the tangential portfolio $\bfw_T$ to exist, the expected return of the minimum variance portfolio $\bfw_{\mvp}$ must be greater than the risk-free
asset return $r_0$,
see \cite{merton1972analytic}. 
This is equivalent to require $(b- r_0 a) >0$.

Let
\[
  \Pi_T := \left(  \begin{array}{c}
  	 \bfw_T \\ 0
  	\end{array} \right)  \qquad \mbox{and} \qquad
  \Pi_0 := \left(  \begin{array}{c}
  	0 \\ 1
  \end{array} \right)
\]
be the respective representation of the tangential portfolio and the risk-free asset in the $(n+1)$ assets space. The portfolios on the CML can then be represented as
\[
  \Pi(\beta) = \Pi_0 + \beta (\Pi_T - \Pi_0) , \quad \beta \ge 0 .
\]
Let $\widetilde{V}$ denote the covariance matrix of the $(n+1)$ assets and
$\widetilde{D}$ be the corresponding Euclidean distance matrix.
We have
\[
  \widetilde{V} = \left[
  \begin{array}{cc}
  	 V & 0 \\
  	 0 & 0
  \end{array} 
  \right] \quad \mbox{and} \quad
  \widetilde{D} =
    \left[
  \begin{array}{cc}
  	D & \frac 12 \eta \\ [0.6ex]
  	\frac 12 \eta^\top & 0
  \end{array} 
  \right] =  \left[
  \begin{array}{cc}
  	\frac 12 (\eta \bfone_n^\top + \bfone_n \eta^\top ) - V & \frac 12 \eta \\ [0.6ex]
  	\frac 12 \eta^\top & 0
  \end{array} 
  \right].
\]
We compute the DR of the portfolio $\Pi$ using $\bfone_n^\top \bfw_T =1$
\begin{eqnarray}
	q(\Pi) &=& \frac 12 \Pi^\top \widetilde{D} \Pi \nonumber \\
	&=& -\frac 12 \beta^2 \bfw_T^\top V \bfw_T 
	  + \frac 12 \beta^2 (\eta^\top \bfw_T) 
	  + \frac 12 \beta(1-\beta) (\eta^\top \bfw_T)  \nonumber \\ 
	&=& -\frac 12 \beta^2 \sigma_T^2 + \frac 12 \beta (\eta^\top \bfw_T) ,
	\label{qPi}
\end{eqnarray}
where $\sigma^2_T$ is the variance of the portfolio $\bfw_T$.
On the other hand, the variance of the portfolio $\Pi$ is 
a quadratic function of
$\beta$: 
\[
  \sigma^2_{\Pi} = \beta^2 \sigma^2_T  .
\]
Substituting it into \eqref{qPi}, we get
\[
  q(\Pi) = -\frac 12 \sigma^2_{\Pi} + \frac{\eta^\top \bfw_T}{2 \sigma_T} \sigma_{\Pi}.
\]
To indicate its dependence on CML, we denote $q(\Pi)$ by 
$q_{\cml}(\sigma)$.
%as the portfolio $\Pi$ is uniquely defined by its standard deviation $\sigma$. 
Therefore, we have
\be \label{q_CML}
  q_{\cml} (\sigma) = -\frac 12 \sigma^2 + \frac{\eta^\top \bfw_T}{2 \sigma_T} \sigma
  = -\frac 12 
  \left(
  \sigma - \frac{\eta^\top \bfw_T}{2 \sigma_T }
  \right)^2 + \frac 18 \left( \frac{\eta^\top \bfw_T}{\sigma_T} \right)^2 .
\ee 
When $\eta^\top \bfw_T >0$, the largest DR portfolio happens at
\[
  \sigma = \frac{\eta^\top \bfw_T }{ 2 \sigma_T} \qquad
  \mbox{or equivalently at} \qquad
  \beta = \frac{ \eta^\top \bfw_T }{2 \sigma^2_T }.
\]

%%%%%% This part is not right AS we cannot compare q_cml with q_dr %%%%%%
%%%%%% Because one involves (n+1) assets and q_dr only involves n assets %%%
%Comparing (\ref{q_CML}) with \eqref{q_efdr} at $\sigma_{\mvp}$, 
%after some simplification we get
%from $q_{\efdr} (\sigma) \ge q_{\cml}(\sigma)$ that 
%\[
% \eta^\top \bfw_T \le (\bfone_n^\top V^{-1} \eta) \sigma_T \sigma_{\mvp}.
%\]
%This inequality suggests a very peculiar scenario that when 
%$\bfone_n^\top V^{-1} \eta$ is negative, then the average risk 
%$\eta^\top \bfw_T$ of the portfolio
%$\bfw_T$ must be negative too, independent of the riskfree return $r_0$
%which determines $\bfw_T$.

Having computed the DR curve \eqref{q_CML} of the CML, we now compute the 
efficient DR curve for the $(n+1)$ assets. 
We note that we cannot directly apply the formula \eqref{q_efdr}) here for the following reasons:
(i) there are $(n+1)$ assets here while  $q_{\efdr}$ has only $n$ risky assets, and 
the DR of a portfolio depends on how many assets it has, see the discussion in Subsection~\ref{Subsection-size};
and (ii) the covariance matrix $V$ for $q_{\efdr}$ is assumed nonsingular, while $\widetilde{V}$ is singular. 
Hence, we need to compute the efficient DR curve from a scratch.
Define
\be \label{wnplus}
  \widetilde{\bfw}(\sigma) 
  := \arg\max \; \widetilde{q}(\widetilde{\bfw}) =
  \frac 12 \widetilde{\bfw}^\top \widetilde{D} \widetilde{\bfw},
  \ \mbox{s.t.} \
  \left\{ 
  \begin{array}{l}
  \bfone_{n+1}^\top \widetilde{\bfw} = 1\\
  \bfw^\top V \bfw = \sigma^2 
  \end{array}  \right.
  \quad \mbox{with} \quad
  \widetilde{\bfw} := \left(
  \begin{array}{c}
  	 \bfw \\ w_0
  \end{array} 
  \right) \in \Re^{n+1}
\ee 
where
$\sigma^2$ is a given level of the risk and
 $w_0$ represents the weight invested in the risk-free asset.

Using the fact that $\bfone_n^\top \bfw = 1 - w_0$, we get
$
   \widetilde{q}(\widetilde{\bfw}) 
   = q(\bfw).
$
Therefore, Problem \eqref{wnplus} is equivalent to the following problem:
\be \label{what}
 \widehat{\bfw}(\sigma) := \arg\max \; q(\bfw) , \quad \mbox{s.t.}
 \quad \bfw^\top V \bfw = \sigma^2 .
\ee 
The weight vector $\widehat{\bfw}(\sigma)$ is the risky part of 
$\widetilde{\bfw}(\sigma)$ and $\widehat{w}_0 := 1 - \bfone_n^\top \widehat{\bfw}(\sigma)$ is the weight invested in the risk-free asset.
Hence, 
$ 
\widetilde{q}(\widetilde{\bfw}(\sigma))
= q (\widehat{\bfw}(\sigma) ).
$
Repeating the computational procedure for $q_{\efdr}$ that leads to 
\eqref{q_efdr}, we can compute for $q (\widehat{\bfw}(\sigma) )$, which is given below:
\be \label{qplus}
  q (\widehat{\bfw}(\sigma) ) = -\frac 12 \sigma^2 + 
  \frac {\sigma}2 \sqrt{\eta^\top V^{-1} \eta}    =: \widetilde{q}_{\efdr}(\sigma),
\ee
where the optimal DR is denoted by $\widetilde{q}_{\efdr}(\sigma)$ to indicates
it is the optimal DR under the given level of risk $\sigma^2$. 
This is the efficient DR curve for the $(n+1)$ assets. 
Comparison of \eqref{qplus} with \eqref{q_CML} leads to the following two remarks.

\begin{itemize}
	\item[(i)] Both curves $\widetilde{q}_{\efdr}(\sigma)$ and $q_{\cml}(\sigma)$ are of standard parabolas, and are much simpler than their counterparts $q_{\efdr}(\sigma)$ and $q_{\ef}(\sigma)$ where
	only risky assets are considered. 
	The fact that  
	$
	 \widetilde{q}_{\efdr}(\sigma) \ge q_{\cml}(\sigma)
	$
	(the efficient DR curve dominates the DR curve of CML)
	yields
	\be \label{New-Inequality}
	  \sqrt{\eta^\top V^{-1} \eta} \sigma_T \ge \eta^\top \bfw_T .
	\ee
	We further note that the quantity $(\eta^\top V^{-1} \eta )$ is $1/\sigma^2_{\eta}$ defined by
	\[
	  \sigma^2_{\eta} := \min \bfw^\top V \bfw, \qquad  \mbox{s.t.} \quad
	  \eta^\top \bfw = 1.
	\] 
	Hence,the inequality \eqref{New-Inequality} becomes
	\[
	  \frac{\sigma_T}{\sigma_{\eta}} \ge \eta^\top \bfw_T .
	\]
	This inequality on the tangential portfolio $\bfw_T$ is new.
	
	\item[(ii)] When the variance vector $\eta$ is proportional to the excess rate of return $(\overline{\bfr} - r_0 \bfone_n)$, i.e.,
	$
	   \eta = \gamma (\overline{\bfr} - r_0 \bfone_n)
	$
	for some $\gamma >0$, then we can prove $q_{\cml}(\sigma) = \widetilde{q}(\sigma)$. In other words,
	the DR curve of CML becomes the efficient DR portfolio. 
	This result generalizes Thm.~\ref{Thm-Two-EF}(iv) to the case a risk-free
	asset is included.
	Otherwise, the gap can be measured as follows:
	\[
	  \widetilde{q}(\sigma) - q_{\cml}(\sigma)
	  = \frac{ \sigma}2 \Big(  \sqrt{\eta^\top V^{-1} \eta} - \frac{\eta^\top \bfw_T }{\sigma_T} \Big) \qquad \forall \ \sigma \ge 0.
	\]
	The gap is a linear function of $\sigma$.
	The curves of $q_{\cml}(\sigma)$ and $\widetilde{q}(\sigma)$ are similar 
	to their counterparts in Fig.~\ref{fig:DR-Frontier}.
	
\end{itemize}

%%%%%%%%%%%%%%%%%%%%%%%%%%%%%%%%%%%%%%%%%%%%%%%%%%%%%%%%%%%%%
\section{A Numerical Illustration} \label{Section-Numerical} 

Theoretically, we are largely clear what to expect of the diversification return based portfolios.
In this part, we illustrate the behaviour of those portfolios using a real data set and 
compare them with some existing portfolios. 
The illustration reveals some interesting observations and
enhances our understanding of those portfolios. 
In particular, we demonstrate the following:

\begin{itemize}
	\item[(i)] We illustrate the concerned portfolios in three graphs: 
	$(\sigma, q)$ (standard deviation and diversification return) graph;
	$(\sigma, c)$ (standard deviation and centrality) graph; and
	$(\sigma, R)$ (standard deviation and return) graph.
	
	\item[(ii)] We also single out some particular portfolios. They include:
	$\mbox{MVP}$: Minimum Variance Portfolio;
	$\MDRP$: Maximum Diversification Return Portfolio;
%	$\mbox{maxVP}$: maximum Volatility Portfolio studied in \cite{maeso2020maximizing}:
%	\[
%	  \bfw_{\maxvp} := \frac{V^{-1}\eta}{\bfone_n^\top V^{-1} \eta }.
%	\]
    $Q$-portfolio: the efficient mean-variance portfolio that has the largest diversification return; and
    $\MDP$: Maximum Diversification ratio Portfolio, studied in \cite{choueifaty2008toward}. This 
    portfolio is explained below.

\end{itemize}

\subsection{Maximum diversification ratio portfolio}

Researchers often get confused between the Maximum Diversification Return Portfolio (MDRP) and 
the Maximum Diversification ratio Portfolio (MDP), which was initially studied
by Choueifaty and Coignard \cite{choueifaty2008toward}:
\be \label{MDP}
 \bfw_{\mdp} := \arg\max \; \frac{\sqrt{\eta }^\top \bfw}{ \sqrt{\bfw^\top V \bfw }  },
 \quad \mbox{s.t.} \quad \bfone_n^\top \bfw = 1.
\ee 
It would be interesting to study their relationship. 
Obviously, \eqref{MDP} is equivalent to
\[
  \bfw_{\mdp} = \arg \max \; \left(\frac{\sqrt{\eta }^\top \bfw}{ \sqrt{\bfw^\top V \bfw }  } \right)^2
  = \frac{ \bfw^\top \sqrt{\eta} \sqrt{\eta}^\top \bfw }{ \bfw^\top V \bfw }
  \quad \mbox{s.t.} \quad \bfone_n^\top \bfw = 1.
\]
Let us fixed the risk level at $\bfw^\top V \bfw = \sigma^2$ and consider the
corresponding maximum diversification ratio
portfolio:
\be \label{mdp-sigma}
  \bfw_{\mdp}(\sigma) := \arg\max\; \frac 12 \bfw^\top \sqrt{\eta} \sqrt{\eta}^\top \bfw,
  \quad \mbox{s.t.} \quad \bfone_n^\top \bfw = 1, \ \
  \bfw^\top V \bfw = \sigma^2. 
\ee 
We estimate the difference of the objectives
of the two portfolios MDRP (\ref{DRP-Risk-2}) and MDP (\ref{mdp-sigma}):
\begin{eqnarray*}
  d_{\sigma} &:=& \frac 12 \max_{\bfw} \; \eta^\top \bfw -
  \frac 12 \max_{\bfw} \; \bfw^\top \sqrt{\eta} \sqrt{\eta}^\top \bfw \\
  &\le&\frac 12  \max_{\bfw} \; \Big( \eta^\top \bfw - \bfw^\top \sqrt{\eta} \sqrt{\eta}^\top \bfw \Big) \\
  &=& \frac 12 \max_{\bfw} 
   \bfw^\top \underbrace{\Big( \frac 12 (\eta^\top \bfone_n + \bfone_n \eta^\top) - \sqrt{\eta} \sqrt{\eta}^\top  \Big)}_{=: D_{\eta}} \bfw 
\end{eqnarray*} 
It is easy to see that
\[
  \diag(D_{\eta}) =0 \quad \mbox{and} \quad -J D_{\eta} J = J \sqrt{\eta} \sqrt{\eta}^\top J \succeq 0,
\]
where $J := I_n - \frac 1n \bfone_n \bfone_n^\top$ is the centering matrix. 
It follows from \cite{gower1985properties} (see also \cite[Eq.(1)]{qi2013semismooth}) that
$D_{\eta}$ is Euclidean distance matrix.
Let us restrict to the long-only portfolio $\bfw \ge 0$ and define
\[
  S_{\sigma}^+ := \left\{
    \bfw \; | \; \bfone_n^\top \bfw = 1, \ \ 
    \bfw^\top V \bfw = \sigma^2, \ \ \bfw \ge 0
  \right\}.
\]
On one hand,
\[
 - 2d_{\sigma} = \max_{\bfw \in  S_{\sigma}^+} \; \bfw^\top \sqrt{\eta} \sqrt{\eta}^\top \bfw
 - \max_{\bfw \in  S_{\sigma}^+} \; \eta^\top \bfw 
 \le \max_{\bfw\in  S_{\sigma}^+} 
 \bfw^\top  (-D_{\eta}) \bfw \le 0,
\]
where the last inequality used the facts that $D_{\eta}$ is EDM (hence its elements are non-negative) and
$\bfw \ge 0$. On the other hand,
\be \label{dmax-Problem}
 d_{\sigma} \le \left\{ \max_{\bfw \ge 0} 
 \frac 12 \bfw^\top  D_{\eta}  \bfw, \quad \mbox{s.t.} \quad \bfone_n^\top \bfw = 1 \right\}
 =: d_{\max}. 
\ee
%%%%%%%%%%%%%%%%%%%%%%%%%%%%%%%%%%%%%%%%%%%%%%%%%%%%%
%We note that $d_{\max}$ is from a quadratic programming (QP) and we implicitly assumed the existence of
%an optimal solution (hence $\max$ was used in the QP).
%We now prove the existence.
%
%%%%%%%%%%%%%%%%%%%%%%%
%\begin{proposition} %%% we cannot prove the result
%Suppose the variance vector $\eta$ is not proportional to $\bfone_n$.
%Then the QP problem (\ref{dmax-Problem}) has a unique optimal solution and $d_{\max}$ is well defined.
%\end{proposition}
%
%{\bf Proof.}
%According to \cite[Thm.~16.2]{nocedal1999numerical}, a sufficient condition for the existence of
%existence of an optimal solution to this QP the nonsingularity of the following
%KKT matrix:
%\[
%  M := \left[
%   \begin{array}{cc}
%   	  D_{\eta} & \bfone_n \\
%   	  \bfone_n^\top & 0
%   	\end{array} 
%  \right]
%\]
%Suppose 
%\[
%  M \left[
%  \begin{array}{c}
%  	 \bfx \\ \lambda
%  \end{array} 
%  \right] = 0  \quad \mbox{or equivalently} \quad
%  \left\{
%   \begin{array}{l}
%    D_{\eta} \bfx + \lambda \bfone_n = 0\\
%    \bfone_n^\top \bfx = 0 
%   \end{array} 
%  \right . \quad \bfx \in \Re^n, \ \ \lambda \in \Re .
%\]
%Let $J= I_n - \frac 1n \bfone_n \bfone_n^\top$ be the projection matrix to the null space of $\bfone_n$.
%Since $\bfx$ is in the null space, the above system is equivalent to
%\[
%  D_{\eta} J \bfy + \lambda \bfone_n = 0.
%\]
%Pre-multiplying by $(J \bfy)^\top$ leads to (using the fact $J bfone_n =0$)
%\[
%   \bfy^\top J  D_{\eta} J \bfy = 0 .
%\]
%%%%%%%%%%%%%%%%%%%%%% proof failed%%%%%%%%%%%%%%%%%%%%%%%%%
Hence, for long-only portfolios (see \cite{qi2021long}), we must have
\[
 \max_{\bfw \in  S_{\sigma}^+} \; \bfw^\top \sqrt{\eta} \sqrt{\eta}^\top \bfw
 \le
 \max_{\bfw \in  S_{\sigma}^+} \; \eta^\top \bfw \le \max_{\bfw \in  S_{\sigma}^+} \; \bfw^\top \sqrt{\eta} \sqrt{\eta}^\top \bfw 
 + d_{\max}.
\]

We note that $d_{\max}>0$ is independent of the level of $\sigma^2$ of the portfolios involved.
This means that, given any risk level $\sigma^2$, the distance between the two portfolios objectives is uniformly bounded
irrelevant to the risk level. 
In other words, the diversification ratio portfolios should follow the trend of the diversification return
portfolios. This is exactly what is observed in the following numerical example.

%%%%%%%%%%%%%%%%%%%%%%%%%%%%%%%%%%%%%%%%%%%%%%%%%%%%%%%%%%%%%%
\subsection{Dax30 data}

\begin{figure}[h!]
	\centering
	\includegraphics[width= 0.65\textwidth]{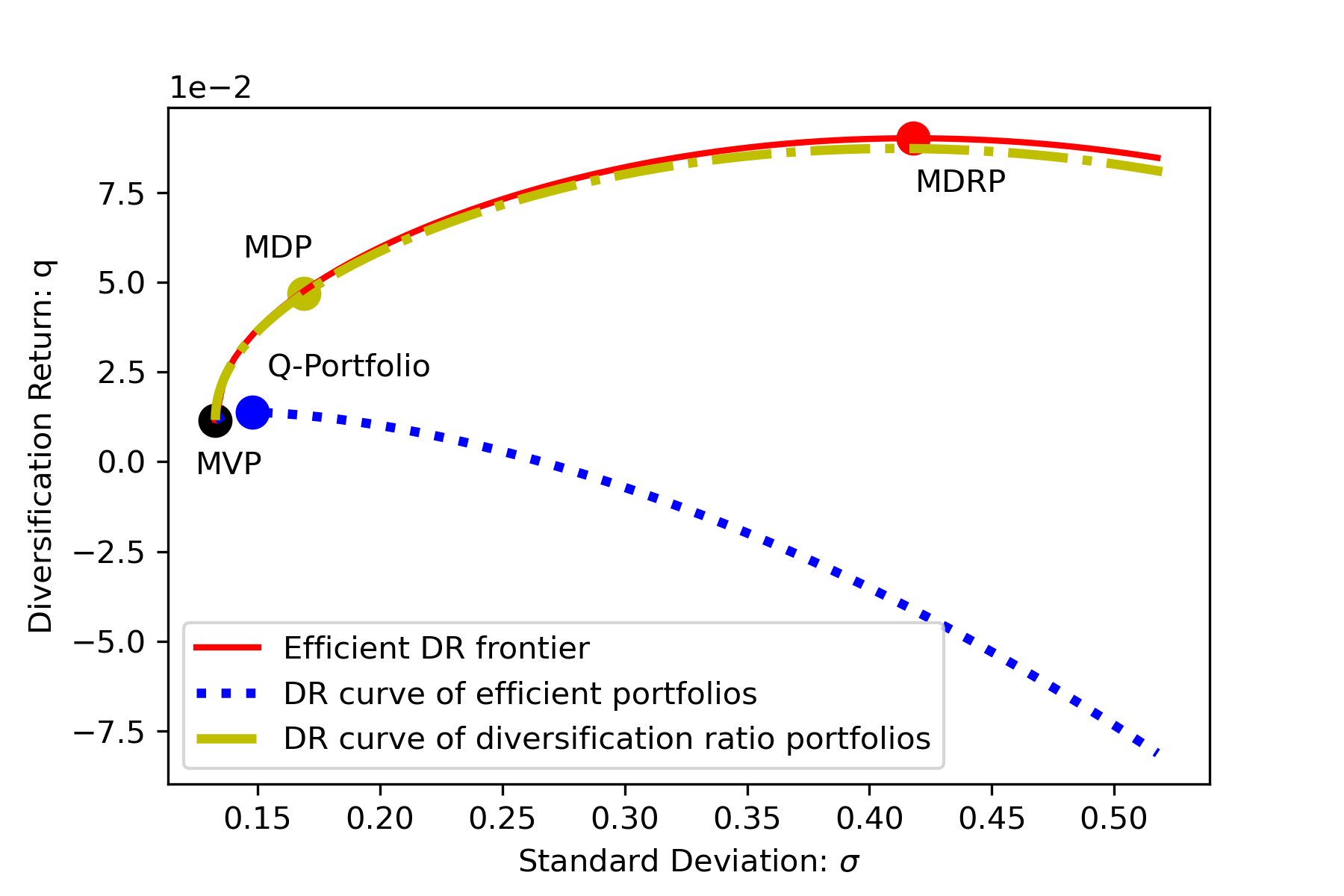}
	\caption{Comparison of portfolios in the $(\sigma, q)$ plane.
	}
	\label{Fig-Dax-DR}
	%\vspace{-16pt}
\end{figure}

\begin{figure}[h!]
	\centering
	\includegraphics[width= 0.65\textwidth]{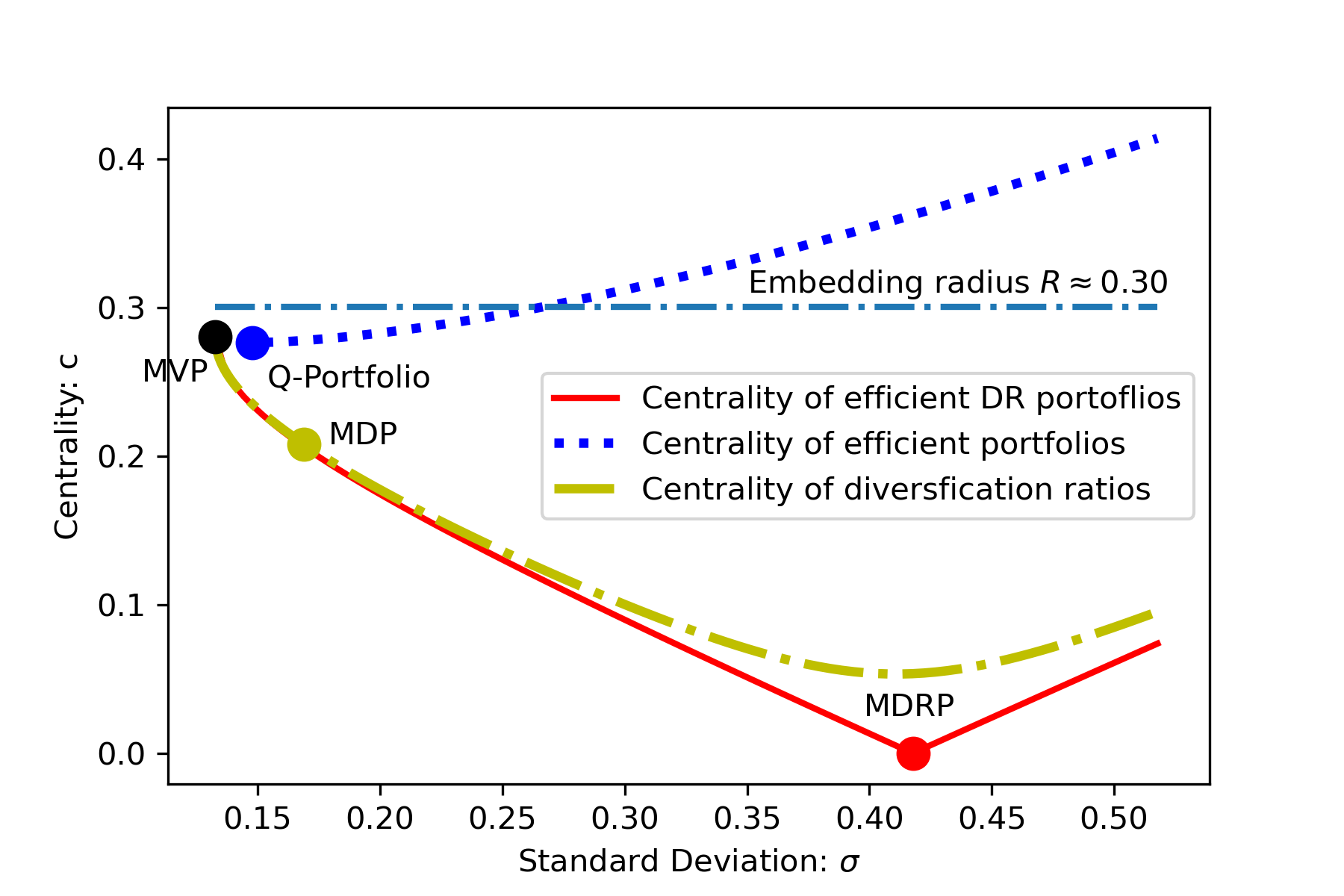}
	\caption{Comparison of portfolios in the $(\sigma, c)$ plane.
	}
	\label{Fig-Dax-Centrality}
	%\vspace{-16pt}
\end{figure}

\begin{figure}[h!]
	\centering
	\includegraphics[width= 0.65\textwidth]{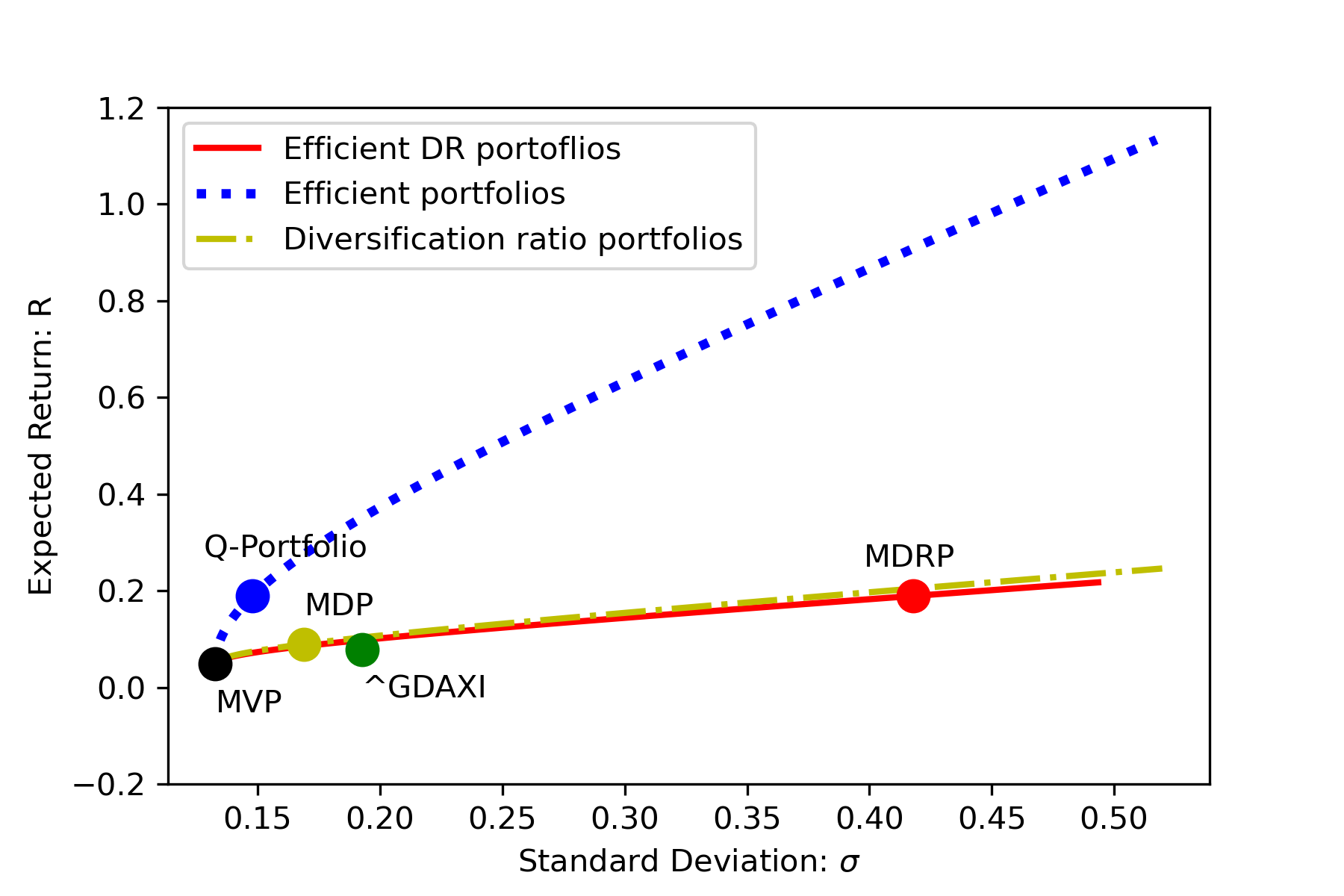}
	\caption{Comparison of portfolios in the $(\sigma, R)$ plane.
	}
	\label{Fig-Dax-Return}
	%\vspace{-16pt}
\end{figure}

This data set consists of $30$ stocks\footnote{The ticker symbols for those
	stocks are
	ADS.DE, ALV.DE, BAS.DE, BAYN.DE, BEI.DE,
	BMW.DE, CBK.DE, CON.DE, DAI.DE, DB1.DE,
	DBK.DE, DPW.DE, DTE.DE, EOAN.DE, FME.DE,
	FRE.DE, HEI.DE, HEN3.DE, IFX.DE, LHA.DE,
	LIN.DE, LXS.DE, 
	MRK.DE, MUV2.DE, RWE.DE,
	SAP.DE, SDF.DE, SIE.DE, TKA.DE, VOW3.DE.} that have appeared in DAX30 Index (DAX30) and
was used in \cite[Page~336]{hilpisch2014python}. 
The data period is from January 3, 2017 to December 31, 2021. 
The mean and the covariance matrix of the daily returns were 
annualized\footnote{Suppose the observations are from time $t_0$ to
	time $T$. Let $n$ denote the number of returns in this period and let $N$ denote the number of calendar days between $t_0$ and $T$. The annualized time step is $\delta = N/(365n)$. Let $\widehat{\mu}$ and $\widehat{V}$ be the
	sample mean and covariance matrix of the returns. Then the annualized mean and 
	covariance matrix are respectively given by $\mu = \widehat{\mu}/\delta$ and
	$V = \widehat{V}/\delta$. The value of $\delta$ for DAX30 data set is
	$0.003952$.}.
The averaged weekly data was also tested and the behaviours of the concerned portfolios are similar
to that of the plotted graphs in Fig.~\ref{Fig-Dax-DR}, Fig.~\ref{Fig-Dax-Centrality}, and Fig.~\ref{Fig-Dax-Return},
 and hence are
omitted.
In all three figures, we plots the portfolios of $\bfw_{\mdrp}(\sigma)$ in (\ref{DRP-Convex}),
$\bfw_{\mv}(\sigma)$ in (\ref{MV-sigma}) and $\bfw_{\mdp}(\sigma)$ in (\ref{mdp-sigma}) as $\sigma$ varies.
They respectively represent the optimal diversification return portfolios, the optimal mean-variance portfolios,
and the diversification return ratio portfolios. 
We note that only $\bfw_{\mv}(\sigma)$ is explicitly related to the stock returns and the other two are only
relevant to the covariances of the stocks.

%%%%%%%%%%%%%%%%%%%%%%%%%%%%%%%%%%%%%%%%%%%%%%%%%%%
%%%%%%%%%%%%%%%%%%%%%%%%%%%%%%%%%%%%%%%%%%%%%%%%%%%%%%%%%%%%%
%%%%%%%%%%%%%%%%%%%%%%%%%%%%%%%%%%%%%%%%%%%%%%%%%%%%%%%%%%%%%%

We summarized the key observations as follows.

\begin{itemize}
\item[(i)] The efficient DR curve of the optimal diversification return portfolio $\bfw_{\mdrp}(\sigma)$ and that of
the optimal diversification ratio portfolio $\bfw_{\mdp}(\sigma)$ are surprisingly close, 
see the standard-deviation and diversification-return  $(\sigma, q)$ graph in
Fig.~\ref{Fig-Dax-DR}. 
It is partially because their respective objectives are not far from each other and the gap is uniformly bounded
irrelevant of the risk level involved. 
Equally surprising is their similar expected returns, as shown in the standard-deviation and return $(\sigma, R)$ graph in
Fig.~\ref{Fig-Dax-Return}. 
However, their difference was vivid in the standard-deviation and centrality $(\sigma, c)$ graph
in Fig.~\ref{Fig-Dax-Centrality}.
It seems that the centrality curve of MDP is smoother than that of MDRP.
Furthermore, the centrality of $\bfw_{\mdp}(\sigma)$ is strictly above that of $\bfw_{\mdrp}(\sigma)$.
However, a general trend is that they follow each other.

\item[(ii)] Contrary to the observation above, the trajectory of the mean-variance portfolio goes the opposite way to
that of MDRP and MDP. 
In particular, the diversification return of efficient MV portfolios quickly went negative, see
Fig.~\ref{Fig-Dax-DR}. 
This means as the standard deviation $\sigma$ (equivalently, the expected return) is above a certain
threshold, they tend to concentrate on few assets leading to less diversified portfolios. 
This is also reflected on the centrality graph Fig.~\ref{Fig-Dax-Centrality}.
Any MV portfolios that are above the horizontal line corresponding to the embedding radius $R \approx 0.30$
has negative $q$ according to the formula \eqref{Centrality-Identity}.
Another drawback of MV portfolios for this dataset is that its expected return grows too fast, see Fig.~\ref{Fig-Dax-Return}.
We also plotted the Dax index return itself, denoted as $^\wedge$GDAXI (the tick symbol of the index), which is very close to the return curve of efficient
diversification return portfolios, and is far below the efficient frontier. 

\item[(iii)] There is a clear cluster among the interested portfolios. 
The portfolios of MVP, MDP and Q-portfolio form a cluster in all three figures. 
Moreover, the market portfolio $^\wedge$GDAXI is also close to this cluster, meaning that portfolios
near this cluster tends to follows the market trend. More numerical experiments are needed to see if
it is a universal observation. Although the efficient portfolios tend to quickly yield negative DR,
the Q-portfolio has positive DR and should be numerically investigated to validate its value.
In contrast, MDRP is far away from the cluster. The return of the efficient DR portfolios in
Fig.~\ref{Fig-Dax-Return} looks flat and this indicates that the increase in return is much slower than that
of the standard deviation. This observation confirms once again that MDRP and MDP are based on two
very different criteria though they both stay close to one curve in each of the plots.

\end{itemize}

%%%%%%%%%%%%%%%%%%%%%%%%%%%%%%%%%%%%%%%%%%%%%%%%%%%%%%%%%%%%%
\section{Conclusion} \label{Section-Conclusion}

This paper provided a thorough study of the diversification return based portfolio from
an optimization perspective.
It shows that there is intrinsic connection between the diversification return and the norm-weighted portfolio, 
and a separation theorem also holds for the efficient diversification return portfolios.
The DR curve of those portfolio is strongly concave in the $(\sigma, q)$ space. Consequently, any
portfolio beyond MDRP is discarded as it would have higher standard deviation than MDRP, but with less
DR. The separation theorem implies that the meaningful portfolios are convex combination of MVP and MDRP.
We also derived a formula for the DR curve of the mean-variance efficient portfolios 
and conducted its comparison with the efficient DR curve. In particular,
the Q-portfolio seems to have some advantages than other MV efficient portfolio.
We also extend such investigation to the case where a risk free asset is also available.
The DR curves for this case are standard parabolas and the risk adjusted gap between the two DR curves
is constant. We also investigated the link between MDRP and the maximum diversification ratio portfolio (MDP). 
Although they are based on different criterion, they tend to follow each other and this is well 
demonstrated in the numerical example. 

The good understanding of the DR portfolios also prompts an important question: how to enhance the DR model under noisy 
environment. This would naturally lead to robust variants, that have been widely studied for the MV models, see
\cite{chen2011tight, ghahtarani2022robust}. We plan to investigate this possibility in our next research.

%%%%%%%%%%%%%%%%%%%%%%%%%%%%%%%%%%%%%%%%%%%%

\bibliographystyle{plain} %{siam}
\bibliography{PortfolioRefs}

\end{document}